\documentclass{article}
\usepackage{graphicx} 
\usepackage{mathtools}
\usepackage{amsfonts}
\usepackage{amsthm}
\usepackage{amsmath}
\usepackage{amssymb}
\usepackage{multirow}
\usepackage{array}
\usepackage{soul}
\usepackage{amssymb}
\usepackage{graphicx}
\usepackage{float}
\usepackage{amsmath}
\usepackage{placeins}
\usepackage{bm}
\usepackage{amsmath}
\usepackage[dvipsnames]{xcolor}
\usepackage{braket}
\usepackage{tikz-cd}
\usepackage[colorlinks,linkcolor=blue]{hyperref}
\usepackage[nameinlink]{cleveref}
\usepackage{stmaryrd}
\usepackage{comment}
\newcommand{\N}{\mathbb{N}}
\newcommand{\R}{\mathbb{R}}
\newcommand{\Q}{\mathbb{Q}}
\newcommand{\Z}{\mathbb{Z}}
\newcommand{\I}{\mathcal{O}}
\newcommand{\C}{\mathbb{C}}

\def\house#1{{%
    \setbox0=\hbox{$#1$}
    \vrule height \dimexpr\ht0+1.4pt width .4pt depth \dp0\relax
    \vrule height \dimexpr\ht0+1.4pt width \dimexpr\wd0+2pt depth \dimexpr-\ht0-1pt\relax
    \llap{$#1$\kern1pt}
    \vrule height \dimexpr\ht0+1.4pt width .4pt depth \dp0\relax
}}

\theoremstyle{definition}
\newtheorem{definition}{Definition}[section]
\newtheorem{example}[definition]{Example}
\newtheorem{conjecture}[definition]{Conjecture}
\theoremstyle{plain}
\newtheorem{theorem}[definition]{Theorem}
\newtheorem{lemma}[definition]{Lemma}
\newtheorem{corollary}[definition]{Corollary}

\newtheorem{proposition}[definition]{Proposition}

\title{On Rankin Invariants of Ideal Lattices}

\author{ Yuan Xiang Chew\thanks{Corresponding author, School of Physical and Mathematical Sciences, Nanyang Technological University, Singapore, email:
\texttt{S240055@e.ntu.edu.sg}.
} \hspace{2mm} and 
Fr\'ed\'erique Oggier
\thanks{School of Mathematics, University of Birmingham, UK.}
}
\date{}

\begin{document}

\maketitle

\begin{abstract}
The Rankin constant $\gamma_{n,k}$ measures the largest discriminant of a densest rank $k$ sublattice of a lattice in $\R^n$ over all such lattices of rank $n$ and their exact values are mostly open. For a general lattice $L\subset \R^n$ of rank $n$, the problem of determining its Rankin invariant $\gamma_{n,k}(L)$ is to find the discriminant of such a densest sublattice. There are no general formulas to compute the discriminant of any arbitrary rank $k$ sublattice. This paper specializes to the case of ideal lattices arising from ideals of number fields, and provides two approaches to study Rankin constants: first, we propose a systematic approach to approximate the minimum discriminant by introducing subideal lattices and proving a closed-form formula to compute $\gamma^{ideal}_{n,k}(L)$, the ideal Rankin invariant of the ideal lattice. Subsequently, we establish necessary conditions on the underlying number field for the existence of an ideal lattice admitting an optimal subideal lattice with a prescribed minimum discriminant. Secondly, we consider ideal lattices with arbitrary sublattices. We improve the current known bounds of the Rankin constants $\gamma_{5,2}$, $\gamma_{7,2}$ and $\gamma_{9,2},\gamma_{9,3},\gamma_{9,4}$, and provide an algebraic construction of the lattice that was used in these dimensions.
\end{abstract}


\section{Introduction}

A \textit{lattice} $L\subseteq \R^n$ of rank $k$ is an integral span of a set  $\set{v_1,\ldots,v_k}\subset L$  of real $n$-dimensional linearly independent basis vectors  where $k\leq n$ and the lattice $L$ is said to be of \textit{full rank} if $k=n$. Let $M$ be a $k\times n$ matrix with basis vectors as rows, then the lattice $L$ is defined as
\[ L = \set{\boldsymbol{\lambda} M | \boldsymbol{\lambda}\in \Z^k}.\] 
We call $M$ a generator matrix of $L$ and  $G = MM^\top$ its corresponding Gram matrix.  The volume and discriminant of $L$ are defined as $\mathrm{vol}(L) = |\sqrt{\det G}|$ and $\mathrm{disc}(L) = \det G$, which do not depend on a choice of a basis. For any $x,y \in \R^n$, let $x\cdot y$ denote the standard dot product in $\R^n$ and let $\|x\|^2= x\cdot x$ be the squared (Euclidean) norm of a vector $x.$ Lastly, we denote the rank of a lattice $L$ as $\mathrm{rank}(L).$ 

\begin{definition}[\protect{\cite[Definition 2.8.3]{Martinet2010PerfectLI}}]\label{def:Rankin}
Let $L\subset \R^n$ be a full rank lattice. For $0 < k < n$, let $L' \subset L$ be any sublattice of rank $k$. Set
\[ 
d_k(L) \coloneqq \min_{L_k \subset L} \operatorname{disc}(L')
\] 
where the minimum is taken over all rank $k$ sublattices $L'\subset L$. 

The {\it Rankin invariant of order $k$} of a lattice $L$ is defined by
\[ 
\gamma_{n,k}(L) \coloneqq \frac{d_k(L)}{\operatorname{disc}(L)^{\frac{k}{n}}}
\] 
and the {\it Rankin constant of order $k$} is 
\[ 
\gamma_{n,k} \coloneqq \sup_{\mathrm{rank}(L)=n} \gamma_{n,k}(L).
\]
\end{definition}

For $k=1$, $\gamma_{n,1}= \gamma_n$ is the Hermite constant of rank $n$, which captures the largest smallest vector of lattices of rank $n$. Values of $\gamma_n$ are known up to dimension 9 and 24 (see \cite{sikiric2025lattice} for the case of dimension 9), and there are known lower and upper bounds, as well as conjectured values of $\gamma_n$ for other dimensions (see \cite{cryptoeprint:2022/677}). Values of $\gamma_{n,k}$ are even less well understood (see. e.g., \cite{oggier2025rankin} for a summary of known and recent results from \cite{sawatani2010note}). Establishing bounds or exact values for $\gamma_{n,k}$ remains difficult. Even for a fixed lattice $L$ of rank $n$, the problem of determining $d_k(L)$ is a difficult problem. It is thus natural to wonder whether restricting the class of lattices considered would help to understand better the density of a sublattice of rank $k$ within a lattice of rank $n$ when the lattice is not arbitrary. 

This paper looks at the class of ideal lattices, which are lattices built from number fields. More specifically, an ideal lattice is the pair $(I,b_\alpha)$ where $I$ is an ideal of a totally real or CM number field $K$ of degree $n$ endowed with a positive definite symmetric bilinear form $b_\alpha$ for some totally positive $\alpha \in K$. Alternatively, the ideal lattice in $\R^n$ can also be realized as the image of $I$ under the so-called canonical embedding. Details can be found in Subsection \ref{subsec:prelim}. The class of ideal lattices is a natural choice to consider as it has been well studied (see e.g., \cite{Bayer-Fluckiger:1999,bayer2002ideal}) and many interesting known lattices have been constructed as ideal lattices (see e.g., \cite{jorge2012rotated,ferrari2022rotated,watanabe2018e8,bayer2005modular}). There are also known constructions of well-rounded ideal lattices (see e.g., \cite{fukshansky2012well,fukshansky2013well,de2019well}) where the norm of the shortest vector is known.

The $k$-Densest Sublattice Problem ($k$-DSP) was introduced in 2006 \cite{gama2006rankin} as a generalization of the Shortest Vector Problem (SVP) which is the case of $k=1$. For a fixed lattice $L \subset \R^n$ of rank $n$, the $k$-DSP is asking for the basis of the sublattice attaining $\gamma_{n,k}(I)$. Much work such as basis reduction algorithms and blockwise lattice reduction have been done since (see e.g.,\cite{dadush2013algorithms,li2014approximating}).

There has also been work on the Rankin Constant from the coding theory view point (see e.g.,\cite{oggier2025rankin}).
\subsection{Contribution and Summary}
\label{subsec:summ}

This paper looks at Rankin invariants of ideal lattices.

In Section \ref{sec: bounds on lattices}, we provide bounds on $d_k(L)$ with a focus on $k=2$ for a general full rank lattice $L\subseteq \R^n.$ Subsequently, we introduce the notion of ideal lattices built from ideals coming from number fields and use the structure of number field to deduce stronger results on $d_2(L)$ where $L$ is an ideal lattice.

Restricting the lattice $L$ to be an ideal lattice, denoted as $(I,b_\alpha)$, does not imply that the sublattices are themselves ideal lattices and for a general sublattice, there are no formula to compute its discriminant. We introduce the notion of subideal lattices in Section \ref{sec:subideal} as a new notion which provides a systematic approach to approximating $d_k(I,b_\alpha)$, and ask whether it is possible to determine the density of a subideal lattice within an ideal lattice through the corresponding ideal Rankin invariant $\gamma^{ideal}_{n,k}(I,b_\alpha)$. We provide a closed-form formula, and look more closely at the case where the lattice is built directly on the ring of integers with $b_\alpha=b_1$, for which we provide a partial answer, as well as applications to lower bound $\gamma_{n,1}$. There is a general bound on $\gamma_{n,1}$ from \cite{cryptoeprint:2022/677} in which our construction gives a sharper lower bound in the case of even dimension between 2 and 24 except 14 and 22.

We next discuss the looseness of $\gamma^{ideal}_{n,k}(I,b_\alpha)$ with respect to $\gamma_{n,k}(I,b_\alpha)$ in Section \ref{sec:loose} and give explicit negative conditions on a number field for equality between both invariants, as well as examples of optimal subideal lattices. We establish that in the case of dimension 6, there are ideal lattice constructions of the $E_6$ lattice for which a subideal lattice of rank 2 attains optimality and for dimension 8, there are ideal lattice constructions of the $E_8$ lattice such that there are subideal lattices of rank 2 and 4 which are optimal. Lastly, we construct an ideal lattice $(I,b_\alpha)$ in dimension 18 that contains an optimal rank 2 subideal lattice attaining a desired $d_2(I,b_\alpha)$ value in Example \ref{Example: zeta27}.

Finally we look at the Rankin invariants for known interesting lattices in Section \ref{sec:lower}. In particular, we look at dimension 5, 7 and 9 and provide an algebraic construction of these lattices and further determine their Rankin invariants. We show that the Rankin invariants for the lattices  $D_5^*$ and $E_7$ improve the known lower bound for $\gamma_{5,2}$ and $\gamma_{7,2}$ found from \cite[Corollary 3]{oggier2025rankin}. For dimension 9, we use the laminated lattice, $\Lambda_9$ which is known to attain $\gamma_{9,1}$ \cite[Corollary 1]{sikiric2025lattice} to provide bounds for $\gamma_{9,2},\gamma_{9,3}$ and $\gamma_{9,4}$.

All computational calculations are performed using SageMath \cite{sage}.

\section{Some Bounds on $\gamma_{n,k}$}\label{sec: bounds on lattices}
\subsection{Approximating $d_2(L)$ for any Lattice $L$}
We will provide some bounds on $\gamma_{n,2}(L)$ where $L\subseteq \R^n$ is a general full rank lattice.
\begin{lemma}\label{Lemma: Lemma for modifying basis}
    Let $L\subseteq \R^n$ be a full rank lattice where $n>2$ and $L'\subseteq L$ be a sublattice of rank 2. There exists a $\Z$-basis $\{x,y\}$ of  $L'$ satisfying 
    \[ \|x\| \leq \|y\| \text{ and } |x\cdot y|\leq \frac{\|x\|^2}{2}.\]
\end{lemma}
\begin{proof}
    Let $(x_0,y_0)$ be an ordered $\Z$-basis of $L'$. Consider $y_1 = y_0-kx_0$ where $k$ is the nearest integer of the ratio $\frac{x_0\cdot y_0}{\|x_0\|^2}$, so $\frac{x_0\cdot y_0}{ \|x_0\|^2} = k + r$ where $|r| \leq \frac{1}{2}$ and $y_1 \in L'$. Then 
    \[ |\frac{x_0\cdot y_0}{\|x_0\|^2} -k|= |r| \leq \frac{1}{2} \iff |x_0\cdot (y_1+kx_0) -k\|x_0\|^2| \leq \frac{\|x_0\|^2}{2}\iff |x_0\cdot y_1|\leq \frac{\|x_0\|^2}{2}.\] Furthermore, $\{x_0,y_1\}$ is still a $\Z$-basis for $L'$ since 
    \[ \begin{pmatrix}
        \leftarrow x_0 \rightarrow \\
        \leftarrow y_1 \rightarrow 
    \end{pmatrix}= \begin{pmatrix}
        1&0\\-k&1
    \end{pmatrix}\begin{pmatrix}
        \leftarrow x_0\rightarrow\\
        \leftarrow y_0 \rightarrow
    \end{pmatrix}\] where $\begin{pmatrix}
        1&0\\
        -k&1
    \end{pmatrix}$ is unimodular.  At this stage, if $\|x_0\| \leq \|y_1\|$, we are done. Otherwise, if $\|x_0\|> \|y_1\|$, we look at the ordered basis $(y_1,x_0)$ and  repeat the argument above to $x_0$ to get $x_1$ such that $|y_1\cdot x_1|\leq \frac{\|y_1\|^2}{2}$ and we know that $\{y_1,x_1\}$ is a $\Z$-basis for $L'$. Again, if $\|y_1\| \leq \|x_1\|$, we are done. Otherwise, $\|y_1\| > \|x_1\|$ and we look at the ordered basis $(x_1,y_1)$ and we apply to argument to $y_1$. We see that
    \[ \|x_0\|> \|y_1\| > \|x_1\|\]where in each iteration, if the process is required to be done again, the norm of the first basis vector is decreasing, and since $L'$ is a lattice, this process will terminate. So in a finite number of iterations, we get a $\Z$-basis $\{x,y\}$ satisfying the condition stated.
\end{proof}
\begin{proposition}\label{Prop: Lattice lower bound for discriminant of rank 2 sublattice}
    Let $L\subseteq \R^n$ be a full rank lattice where $n >2$ and $L'\subset L$ be a sublattice of rank 2. Then
    \[ \operatorname{disc}(L') \geq \frac{3}{4}d_1(L)^2.\]
\end{proposition}
\begin{proof}
Let $x \in L'$ be such that $\|x\|^2= d_1(L')$. Then there exists a $\Z$-basis $\{x,u\}$ \cite{magliveras2008primitive} for $L'$. By Lemma \ref{Lemma: Lemma for modifying basis}, applying the algorithm in the proof once to $u$, there exists $y \in L'$ such that $\{x,y\}$ is a $\Z$-basis for $L'$ where 
\[ \|x\|\leq \|y\| \text{ and } |x\cdot y| \leq \frac{\|x\|^2}{2}= \frac{d_1(L')}{2}.\] Then 
\[ \operatorname{disc}(L') = \|x\|^2\|y\|^2- (x\cdot y)^2 \geq d_1(L')^2- (\frac{d_1(L')}{2})^2\geq \frac{3}{4}d_1(L)^2\] where the last inequality is due to the fact that $d_1(L)\leq d_1(L')$ since $L'\subset L.$
\end{proof}
Proposition \ref{Prop: Lattice lower bound for discriminant of rank 2 sublattice} generalizes \cite[Lemma 2]{oggier2025rankin} where given any integral even lattice $L$ (a lattice $L$ whose Gram matrix has integer coefficients and whose diagonal coefficients are even) with $d_1(L)= 2k$ for some positive integer $k$, if $L'\subset L$ is a sublattice of rank 2, then 
\[ \operatorname{disc}(L') \geq \frac{3}{4}d_1(L)^2 = \frac{3}{4}(2k)^2 =3k^2 \geq 3\] and hence, $d_2(L) \geq 3$. However, as seen in the proof of \cite[Lemma 2]{oggier2025rankin}, it is still useful to know that $d_2(L)\equiv 0,3 \pmod 4 $ to further narrow down the possible values.

\begin{definition}[\protect{\cite[Definition 2.6.7]{Martinet2010PerfectLI}}]
    Let $L\subseteq \R^n$ be a full rank lattice. For $1\leq k \leq n$, define the $k$-th successive minimum of $L$ to be 
    \[ \lambda_k(L) = \inf\{r>0  \mid \exists \{x_1,\ldots,x_k\}\subset L \text{ linearly independent and } \|x_i\|^2\leq  r ~ \forall i\}.\]
\end{definition}
\begin{proposition}\label{Prop: general bounds on d_2(L)}
     Let $L\subseteq \R^n$ be a full rank lattice where $n >2$. Then
     \[\frac{3d_1(L)^2}{4}\leq d_2(L)\leq \lambda_2(L)^2.\] 
\end{proposition}
\begin{proof}
   By the definition of the 2nd successive minimum, there exists $\{x,y\}\subset L$ linearly independent and both $\|x\|^2,\|y\|^2 \leq \lambda_2(L).$ Let $L'$ be the sublattice spanned by $\{x,y\}$ of rank 2, then 
   \[ d_2(L) \leq \mathrm{disc}(L')=\|x\|^2\|y\|^2 - (x\cdot y)^2\leq \|x\|^2\|y\|^2 \leq \lambda_2(L)^2.\] For the lower bound, by Proposition \ref{Prop: Lattice lower bound for discriminant of rank 2 sublattice}, $\frac{3}{4}d_1(L)^2$ is a uniform lower bound for any sublattice of rank 2, hence, 
   \[ d_2(L)= \min_{\substack{L'\subset L\\ \mathrm{rank}(L') =2}} \mathrm{disc}(L') \geq \frac{3d_1(L)^2}{4}.\]
 \end{proof}

 \begin{corollary} \label{Cor: Possible values of d_2(L)}
 Let $L\subseteq \R^n$ be a full rank lattice where $n >2$. If there exists a linearly independent pair of vectors with squared norm $d_1(L)$, then
     \[\frac{3d_1(L)^2}{4}\leq d_2(L)\leq d_1(L)^2.\] 
 \end{corollary}
 \begin{proof}
Since there exists a linearly independent pair of vectors with squared norm $d_1(L)=\lambda_1(L)$, then $\lambda_2(L) = \lambda_1(L)$ and by Proposition \ref{Prop: general bounds on d_2(L)}, 
\[ \frac{3d_1(L)^2}{4}\leq d_2(L)\leq \lambda_1(L)^2= d_1(L)^2.\]
 \end{proof}

 The bounds given in Corollary \ref{Cor: Possible values of d_2(L)} are tight. For instance, the lattices $\Z^n$ and $D_n$ for $n>2$ both have a linearly independent set $\{x,y\}$ of lattice points attaining $d_1(\Z^n)= 1$ and $d_1(D_n)= 2$ respectively. Then $d_2(\Z^n)= 1$ attained by the sublattice $\Z^2$ where $d_2(\Z^n)= d_1(\Z^n)^2 =1 ,$ the upper bound. On the other hand, $d_2(D_n)= 3 = \frac{3}{4}d_1(D_n)^2$, shown in \cite[Proposition 6]{oggier2025rankin}, the lower bound. 
 
 In the case of an integral even lattice $L$ containing a linearly independent pair of lattice points with squared norm $d_1(L)$, by Corollary \ref{Cor: Possible values of d_2(L)} and \cite[Lemma 2]{oggier2025rankin}, the possible values of $d_2(L)$ are restricted. Some possible values when $d_1(L)$ is small are given in Table \ref{tab:d1small}.
 
 \begin{table}[h]
\centering
\begin{tabular}{cc}
$d_1(L)$& possible values of $d_2(L)$ \\
\hline
2 & $3,\,4$ \\
4 &  $12,\,15,\,16$ \\
6 &$27,\,28,\,31,\,32,\,35,\,36$ \\
8 &  $48,\,51,\,52,\,55,\,56,\,59,\,60,\,63,\,64$ \\
\end{tabular}
\caption{\label{tab:d1small}
Values of $d_1(L),d_2(L)$ for an integral even lattice $L$ containing a linearly independent pair of points with squared norm $d_1(L)$.
}
\end{table}

Specifically in the case where $k=2,$ for a fixed dimension $n>2$, the Rankin invariant of any full rank lattice $L$ will provide a lower bound for $\gamma_{n,2}$. We would like to have a benchmark to know whether the choice of a lattice $L$ gives a good approximation. In the following result, we will show that it is possible to use a known lattice of rank $\ell< n$ to approximate $\gamma_{n,2}.$
\begin{proposition}\label{Prop: Approximating using lower rank lattice}
    Let $n >2$ and $2\leq \ell < n$ such that $L\subseteq \R^{\ell}$ is a full rank lattice that contains at least two linearly independent lattice points of squared norm $d_1(L)$. Then the lattice $\Lambda= L \oplus \sqrt{d_1(L)}\Z^{n-\ell} \subseteq \R^n$ is a full rank lattice that satisfies 
    \[\frac{3}{4}\gamma_{\ell,1}(L)^\frac{2\ell}{n}\leq \gamma_{n,2}(\Lambda)\leq \gamma_{\ell,1}(L)^\frac{2\ell}{n}.\] In particular,
    \[\gamma_{n,2}\geq\frac{3}{4}\gamma_{\ell,1}^\frac{2\ell}{n}.\]
\end{proposition}
\begin{proof}
    Consider the lattice $\Lambda = L \oplus \sqrt{d_1(L)}\Z^{n-\ell}$ where \[\mathrm{disc}(\Lambda)= \mathrm{disc}(L)\mathrm{disc}(\sqrt{d_1(L)}\Z^{n-\ell})= d_1(L)^{n-\ell}\mathrm{disc}(L).\] Firstly, we note that by Corollary \ref{Cor: Possible values of d_2(L)}, $d_2(\Lambda)\geq \frac{3}{4}d_1(\Lambda)^2$. Hence this gives
    \[ \gamma_{n,2}(\Lambda) =\frac{d_2(\Lambda)}{\mathrm{disc}(\Lambda)^{2/n} }\geq  \frac{3}{4}(\frac{d_1(\Lambda)^2}{(d_1(L)^{n-\ell}\mathrm{disc}(L))^\frac{2}{n}}).\] We prove that $d_1(\Lambda)=d_1(L)$. To see this, let $x \in \Lambda\backslash\{0\}$ where we write $x=(x_1,x_2)$ such that $x_1 \in L$ and $x_2 \in \sqrt{d_1(L)}\Z^{n-\ell}$. Then 
    \[ \|x\|^2 = (x_1\cdot x_1) + (x_2\cdot x_2) = \|x_1\|^2 + \|x_2\|^2.\] Since $\|x_1\|^2, \|x_2\|^2 \geq d_1(L)$, then $d_1(\Lambda)$ is attained when either $\|x_1\|^2 =d_1(L)$ and $\|x_2\|^2 =0$ or $\|x_1\|^2 =0$ and $\|x_2\|^2 = d_1(L)$. Hence, $d_1(\Lambda)= d_1(L)$ and combining with the fact that $\gamma_{\ell,1}(L)= \frac{d_1(L)}{\mathrm{disc}(L)^\frac{1}{\ell}}$, then 
    \begin{equation*}
        \begin{split}
            \gamma_{n,2}(\Lambda) &\geq \frac{3}{4}(\frac{d_1(L)^2}{(d_1(L)^{n-\ell}\mathrm{disc}(L))^\frac{2}{n}})\\
            &=\frac{3}{4}(\frac{d_1(L)^{\frac{2\ell}{n}}}{\mathrm{disc}(L)^\frac{2}{n}})\\
            &= \frac{3}{4}(\frac{d_1(L)}{\mathrm{disc}(L)^\frac{1}{\ell}})^\frac{2\ell}{n}\\
            &= \frac{3}{4}\gamma_{\ell,1}(L)^\frac{2\ell}{n}.
        \end{split}
    \end{equation*} For the upper bound, we see that $d_2(\Lambda)\leq d_1(L)^2$ as $L$ contains at least two linearly independent lattice points of squared norm $d_1(L)$, hence 
    \[ \gamma_{n,2}(\Lambda)\leq \frac{d_1(L)^2}{(d_1(L)^{n-\ell}\mathrm{disc}(L))^\frac{2}{n}} = \gamma_{\ell,1}(L)^\frac{2\ell}{n}.\] In particular, we have that 
    \[ \gamma_{n,2}\geq \sup_{\substack{\mathrm{rank}(L)=\ell}} \frac{3}{4}\gamma_{\ell,1}(L)^\frac{2\ell}{n} =\frac{3}{4}\gamma_{\ell,1}^\frac{2\ell}{n}.\]
\end{proof}

\begin{example}
    For instance, for $n=5$, since $\gamma_{4,1} = \sqrt{2}$, then
\[ \gamma_{5,2}\geq \frac{3}{4}\gamma_{4,1}^\frac{8}{5}\approx 1.3058.\] This serves as a starting point for the search of rank 5 lattices achieving $\gamma_{5,2}$. Using the lattice $D_5$, we have $\gamma_{5,2}(D_5)= \frac{3}{4^\frac{2}{5}}\approx 1.7230$ \cite[Proposition 6]{oggier2025rankin} which provides a new benchmark. In Section \ref{sec:lower}, we will show that there is a lattice that outperforms $D_5.$
\end{example}

\begin{example}\label{Ex: dim 9 and 25}
    Consider $n=9$ and $25$. Then we have that 
    \[ \gamma_{9,2}\geq \frac{3}{4}\gamma_{8,1}^\frac{16}{9}\approx 2.5717 \text{ and } \gamma_{25,2}\geq \frac{3}{4}\gamma_{24,1}^\frac{48}{25}\approx 10.7403.\] This serves as a benchmark for the search of lattices of rank $9$ and $25$ that outperform these lower bounds.
\end{example}

However, from Proposition \ref{Prop: Approximating using lower rank lattice}, we see that the approximation becomes poorer as $n$ increases. This suggests that we should use a lattice $L$ of rank $\ell = n-1$ with large Rankin invariant $\gamma_{\ell,1}(L)$ to provide a lower bound for $\gamma_{n,2}$.

\subsection{Ideal Lattices}

\label{subsec:prelim}

For a number field $K$ of degree $n$, denote by $\mathrm{Emb}_\Q(K,\C)=\set{\sigma_1,\ldots, \sigma_n}$ the set of embeddings of $K$ containing field homomorphisms of the form $\sigma_i: K \rightarrow \C$ fixing $\Q$, $1\leq i \leq n$.
The field trace and norm are respectively given by $\mathrm{Tr}_{K/\Q}(x) = \sum_{i=1}^n \sigma_i(x)$ and  $\mathrm{N}_{K/\Q}(x)= \prod_{i=1}^n \sigma_i(x)$, for any $x \in K$ where $\sigma_i \in \mathrm{Emb}_\Q(K,\C).$
An embedding $\sigma_i \in \mathrm{Emb}_\Q(K,\C)$ is called a real embedding if $\sigma_i(K) \subseteq \R$ and is called a complex embedding if $\sigma_i(K) \not\subseteq \R.$ 
A number field $K$ is called totally real if $\mathrm{Emb}_\Q(K,\C)$ only contains real embeddings.
A number field $K$ is called CM if there exists a totally real number field $F \subseteq K$ where $[K:F]= 2$ and $\sigma_i(K)\not\subseteq \R$ for all $i=1,\ldots,n.$ For any $x\in K\subset \C$, let $\overline{x}$ denote the complex conjugate of $x.$ We say that an element $x \in K$ is totally positive if $\sigma_i(x)>0$ for all $i=1,\ldots,n$.

\begin{definition}[\cite{bayer2002ideal}]\label{Def: Ideal lattice from number field}
Let $K$ be a totally real or CM number field of degree $n$. An {\it ideal lattice} is a pair $(I,b_{\alpha})$ where $I$ is a fractional ideal of $K$, $\alpha \in K$ is totally positive and $b_\alpha$ is a positive definite symmetric bilinear form given by
\begin{equation*}
    \begin{split}
        b_\alpha: I \times I &\rightarrow  \R \\
        (x,y) &\mapsto \mathrm{Tr}_{K/\Q}(\alpha x \overline{y}).
    \end{split}
\end{equation*}
\end{definition}

By the above definition, for any $\lambda \in K$, then $b_\alpha(\lambda x,y)= b_\alpha(x,\overline{\lambda}y)$.

Set $I_\R= I\otimes_\Z \R$ and $K_\R = K\otimes_\Q \R$. Since $I$ is a fractional ideal of $K$, then $I$ contains a $\Q$-basis for $K$ and hence, 
\[ I_\R = I\otimes_\Z \Q \otimes_\Q \R = K_\R.\] We can extend $b_\alpha$ to $K_\R$ by letting 
\begin{equation*}
    \begin{split}
        b_\alpha: K_\R \times K_\R &\rightarrow \R\\
         (x\otimes r_1, x\otimes r_2) &\mapsto r_1r_2 \mathrm{Tr}_{K/\Q}(\alpha x\overline{y})
    \end{split}
\end{equation*}  which is a positive definite symmetric $\R$-bilinear form when $K$ is totally real or CM and $\alpha\in K$ is totally positive. In this paper, the only time we will view $b_\alpha$ as a bilinear form on $K_\R$ is in Subsection \ref{subsec: dimension 5}, where we consider the notion of dual lattice.

A generator matrix for $(I,b_\alpha)$ is obtained via the so-called canonical embedding $\sigma_\alpha$ \cite[p.82]{Bayer-Fluckiger:1999} of a basis of $I$. 
Set  $\alpha_i = \sqrt{\sigma_i(\alpha)}$. 
If $K$ is totally real, this canonical embedding is
\begin{equation*}
    \begin{split}
        \sigma_\alpha: K &\rightarrow \R^n\\
        x &\mapsto (\alpha_1\sigma_1(x),\ldots,\alpha_n\sigma_n(x))
    \end{split}
\end{equation*} and if $K$ is CM, we either use the complex canonical embedding
\begin{equation*}
    \begin{split}
        \sigma_\alpha: K &\rightarrow \C^n\\
        x &\mapsto (\alpha_1\sigma_1(x),\ldots,\alpha_n\sigma_n(x))
    \end{split}
\end{equation*} 
or order the embeddings such that $\sigma_{j}$ and $\sigma_{n/2+j}$ are pairs of complex embeddings for each $j \in \set{1,\ldots,n/2}$ and define the real canonical embedding 
\begin{equation*}
    \begin{split}
    \sigma_\alpha: K &\rightarrow \R^{n}\\
    x&\mapsto \sqrt{2}(\alpha_1\Re\sigma_1(x),\alpha_1 \Im \sigma_1(x),\ldots ,\alpha_{n/2}\Re\sigma_{n/2}(x),\alpha_{n/2} \Im \sigma_{n/2}(x)).
    \end{split}
\end{equation*} These two real canonical embeddings satisfy $\sigma_\alpha(x)\cdot \sigma_\alpha(y)= \mathrm{Tr}_{K/\Q}(\alpha x\overline{y})$ for any $x,y \in K.$ Consequently, we may equivalently write $\sigma_\alpha(I)$ or $(I,b_\alpha)$ which is a full rank lattice in $\R^n$. In fact, the fractional ideal structure in $K$ is not required to construct a lattice. Given any free $\Z$-submodule $M\subset K$ of rank $n$, $\sigma_\alpha(M)$ is a full rank lattice in $\R^{n}$. In this paper, we distinguish the following two notations for a lattice arising from an ideal or a free $\Z$-module that is not an ideal: we write $(M,b_\alpha)$ when $M$ is an ideal of $K$, and $\sigma_\alpha(M)$ when $M$ is a free $\Z$-module but not an ideal.

For an integral ideal $I$, that is $I\subset \I_K$, its norm is defined to be
    \[ \mathrm{N}_{K/\Q}(I) = [\I_K:I]\] which is finite.
In the case where $I$ is a fractional ideal, its norm is  
    \[ \mathrm{N}_{K/\Q}(I) = \frac{\mathrm{N}_{K/\Q}(cI)}{\mathrm{N}_{K/\Q}(c)}\] for some $c \in \I_K$ such that $c I \subset \I_K$.

\begin{proposition}[\protect{\cite[Proposition 2]{bayer2002ideal}}]
    For an ideal lattice $(I,b_\alpha)$ of $K$, its discriminant is given by 
    \[\operatorname{disc}(I,b_\alpha) = \mathrm{N}_{K/\Q}(\alpha)\mathrm{N}_{K/\Q}(I)^2|\Delta_K|\] where $\Delta_K$ is the discriminant of the number field $K$.
\end{proposition} 
\begin{definition}
   Let $(I,b_\alpha)$ be an ideal lattice of $K$. The {\it Rankin invariant of order $1\leq k < n$} of an ideal lattice $(I,b_\alpha)$ is 
    \[ \gamma_{n,k}(I,b_\alpha) =\frac{d_k(I,b_\alpha)}{(\mathrm{N}_{K/\Q}(\alpha)\mathrm{N}_{K/\Q}(I)^2|\Delta_K|)^{\frac{k}{n}}}
    \]
where $d_k(I,b_{\alpha}) = \min_{\substack{L \subseteq (I,b_\alpha)\\ \mathrm{rank}(L)=k}} \operatorname{disc}(L)$.
\end{definition}

For a number field $K$, let $\mu_K\subset \I_K^\times$ denote the subgroup containing all roots of unity of $K.$ It is a well known fact that $\mu_K$ is a finite cyclic group. 
\begin{proposition}\label{Prop: general bounds in CM case}
    Let $K$ be a CM field such that $|\mu_K|>2$. For any ideal lattice $(I,b_\alpha),$ then
    \[ \frac{3d_1(I,b_\alpha)^2}{4} \leq d_2(I,b_\alpha)\leq d_1(I,b_\alpha)^2.\] 
\end{proposition}
\begin{proof}
    It suffices to show that there exists $\{x,y\}\subset I$ that is linearly independent where $b_\alpha(x,x)=b_\alpha(y,y)= d_1(I,b_\alpha)$ and use Corollary \ref{Cor: Possible values of d_2(L)} to attain the bound. Since $|\mu_K|> 2$, there exists $\xi \in \mu_K$ such that $\xi \neq \pm 1$ which implies that $\{1,\xi\}$ is linearly independent. Let $x\in I$ be such that $b_\alpha(x,x)= d_1(I,b_\alpha)$, then $\{x,x\xi\}\subset I$ is also linearly independent, otherwise, there would exist non-trivial coefficients $a,b\in \Z$ such that 
    \[ ax +bx\xi = 0 \iff x^{-1}(ax +bx\xi) =x^{-1}(0)=0 \iff a+b\xi = 0\] which contradicts the fact that $\{1,\xi\}$ is linearly independent. Lastly, 
    \begin{equation*}
        \begin{split}
           b_\alpha(x\xi,x\xi)&= \mathrm{Tr}_{K/\Q}(\alpha (x\xi)\overline{(x\xi)})\\
           &=  \mathrm{Tr}_{K/\Q}(\alpha x\xi \overline{\xi}\overline{x})\\
            &= \mathrm{Tr}_{K/\Q}(\alpha x\overline{x})\\
            &= d_1(I,b_\alpha).
        \end{split}
    \end{equation*}
\end{proof} 
For a totally real or CM field $K$, we note that the complex conjugation is an automorphism of $K$ (the identity map on totally real field $K$). Denote $K^+ =K\cap \R$ to be the maximal real subfield of $K$ ($K^+=K$ if $K$ is totally real).
\begin{lemma}\label{Lemma: totally positive elements are in maximal real subfield}
    Let $K$ be a number field and let $\alpha \in K$ be totally positive. Then $\alpha \in K^+.$
\end{lemma}
\begin{proof}
   A particular embedding of $K$, say $\sigma_1$ is the identity map on $K$ and since $\alpha$ is totally positive, this implies that $\sigma_1(\alpha)=\alpha >0$, and hence, $\alpha \in K \cap \R = K^+$. 
\end{proof}
\begin{proposition}\label{Prop: d_2 when 6th root of unity is present}
    Let $K$ be a CM field and $(I,b_\alpha)$ be an ideal lattice. If there exists $\xi \in \mu_K\backslash\{\pm 1\}$ such that $\xi^6=1$, then 
    \[ d_2(I,b_\alpha)= \frac{3}{4}d_1(I,b_\alpha)^2.\]
\end{proposition}
\begin{proof}
    From Proposition \ref{Prop: general bounds in CM case}, we know that 
    \[ d_2(I,b_\alpha)\geq \frac{3}{4}d_1(I,b_\alpha)^2\] and we will prove that there exists a sublattice $L\subset (I,b_\alpha)$ of rank $2$ such that $\mathrm{disc}(L)= \frac{3}{4}d_1(I,b_\alpha)^2$. Let $x \in I$ be such that $b_\alpha(x,x) = d_1(I,b_\alpha)$ and $\xi \in \mu_K\backslash \{\pm 1\}$ such that $\xi^6 = 1$. We know that $\{x,x\xi\}$ is linearly independent and consider the sublattice $L\subset (I,b_\alpha)$ of rank 2 spanned by the image of $x,x\xi $ under $\sigma_\alpha.$ Set $t = \alpha x\overline{x}$. As $K$ is CM, the complex conjugation is an automorphism of $K$ which means both $x,\overline{x}\in K$, implying that that $x\overline{x}\in K\cap \R = K^+.$ By Lemma \ref{Lemma: totally positive elements are in maximal real subfield}, $\alpha \in K^+$ and so, $\alpha x\overline{x}\in K^+$. We have that
    \begin{equation*}
        \begin{split}
            b_\alpha(x,x\xi) &=  Tr_{K/\Q}(\alpha x \overline{x\xi})\\ 
            &= Tr_{K/\Q}(t\overline{\xi})\\
            &= Tr_{K^+/\Q}(t Tr_{K/K^+}(\xi))\\
            &= Tr_{K^+/\Q}(t(\xi+\xi^{-1})).
        \end{split}
    \end{equation*} Since $\xi $ is a $6$-th root of unity, then it is a root of the $6$-th cyclotomic polynomial $X^2-X+1$, that is $\xi^2-\xi +1 = 0.$ Dividing by $\xi$, we get 
    \[ \xi - 1 +\xi^{-1}= 0 \iff \xi+ \overline{\xi} = 1\] and hence, 
    \[ b_\alpha(x,x\xi)= Tr_{K^+/\Q}(t).\] Lastly, we observe that 
    \[ d_1(I,b_\alpha)=  Tr_{K/\Q}(\alpha x\overline{x}) = 2Tr_{K^+/\Q}(t)\] and so 
    \[ b_\alpha(x,x\xi)= \frac{d_1(I,b_\alpha)}{2}.\] Finally, the discriminant of the sublattice $L$ is
    \[ \mathrm{disc}(L) = b_\alpha(x,x)b_\alpha(x\xi,x\xi) - b_\alpha(x,x\xi)^2 = \frac{3d_1(I,b_\alpha)^2}{4}.\]
\end{proof}

\begin{example}
    Let $K = \Q(\zeta_{42})$ be a cyclotomic field of degree 12 where $\zeta_{42}$ is a primitive $42$-th root of unity and since $6 \vert 42$, there exists a primitive $6$-th root of unity in $K$. Consider the ideal lattice $(\I_K,b_\alpha)$ where $\alpha =1$.  By \cite[Lemma 3.2]{fukshansky2012well} we know $d_1(\I_K,b_\alpha) = 12$. By Proposition \ref{Prop: d_2 when 6th root of unity is present}, we know that 
    \[ d_2(\I_K,b_\alpha) = \frac{3(12)^2}{4}= 108 .\] Hence, 
    \[ \gamma_{12,2}(\I_K,b_\alpha)= \frac{108}{(|\Delta_K|)^\frac{2}{12}}\approx 1.4054.\]
\end{example}

In most cases, the CM fields we will be looking at are cyclotomic fields of degree $n > 2$, hence, satisfying Proposition \ref{Prop: general bounds in CM case} immediately. Subsequently, $\zeta_m$ will always denote a primitive $m$-th root of unity.


\section{The Case of Subideal Lattices}
\label{sec:subideal}

Let $K$ be a totally real or CM number field of degree $n$ and $\alpha \in K$ be totally positive. For an ideal lattice $(I,b_\alpha)$, we note that for a sublattice $L\subseteq (I,b_\alpha)$ of rank $k$, $\sigma^{-1}_\alpha(L)\subseteq I$ is a free $\Z$-module of rank $k$. So, a sublattice of an ideal lattice $(I,b_\alpha)$ of rank $k$ corresponds to a free $\Z$-submodule $M\subseteq I$ of rank $k.$ We define the ratio
\[ \gamma_{n,k}((I,b_\alpha),\sigma_\alpha(M)) = \frac{\operatorname{disc}(\sigma_\alpha(M))}{(\mathrm{N}_{K/\Q}(\alpha)\mathrm{N}_{K/\Q}(I)^2|\Delta_K|)^{\frac{k}{n}}}.\]
We first consider the particular case where $L$ is itself an ideal lattice, that is the sublattice corresponds to a fractional ideal of the number field, which we call a subideal lattice.

\begin{definition}
     Let $(I,b_\alpha)$ be an ideal lattice of $K$ and  $1\leq k <n$. A {\it $k$-subideal} lattice $L$ is a sublattice of $(I,b_\alpha)$ of rank $k$ that corresponds to a fractional ideal in a subfield $F\subseteq K$ of degree $k$, i.e., there exists some fractional ideal $J$ of $F$ such that $L = \sigma_\alpha(J).$ We thus denote it as $L=(J,b_\alpha|_F)$ where $b_\alpha|_F$ is the restriction of the bilinear form $b_\alpha$ to the subfield $F$.
\end{definition}

\begin{definition}\label{def:gammaideal}
For an ideal lattice $(I,b_{\alpha})$ of $K$, the {\it ideal Rankin invariant of order $k$} is defined to be
    \[
    \gamma^{ideal}_{n,k}(I,b_\alpha) = \frac{\min_{L\subseteq (I,b_\alpha)}\operatorname{disc}(L)}{\operatorname{disc}(I,b_\alpha)^{k/n}}
    \]
    where $L=(J,b_{\alpha}|_F)$ is a $k$-subideal lattice of $(I,b_\alpha)$ and $F$ is a subfield of degree $k$ of $K$. If there is no subfield of degree $k$ in $K$, we set $\gamma^{ideal}_{n,k}(I,b_\alpha)= \infty.$
\end{definition}

In Subsection \ref{subsec: closed form formula for ideal rankin invariant}, we will show that $\min_{L_k\subset (I,b_\alpha)}\operatorname{disc}(L_k)$ is well-defined and there exists a $k$-subideal lattice attaining it. Clearly, 
\[ \gamma_{n,k}(I,b_\alpha)\leq \gamma^{ideal}_{n,k}(I,b_\alpha).\]
\subsection{A Closed-Form Formula for $\gamma^{ideal}_{n,k}(I,b_\alpha)$}\label{subsec: closed form formula for ideal rankin invariant}

We first show that given a subfield $F\subset K$ of degree $k$, the set of $k$-subideal lattices of $(I,b_\alpha)$ is not empty.

Let $\I_K$ and $\I_F$ be the rings of integers of $K$, respectively $F$. Since $\I_F$ is a Noetherian ring, every $\I_F$-module $M$ which is finitely generated is Noetherian \cite[Chapter X, Proposition 1.4]{lang2012algebra}, which, by definition \cite[p.413]{lang1994algebraic}, means that every $\I_F$-submodule of $M$ is finitely generated.

\begin{proposition}\label{Lemma: Existence of ideal below}
    Let $F$ be a subfield of $K$ where $[F:\Q]=k < n$ and $I$ be a fractional ideal in $K$. The set of $k$-subideal lattices of $(I,b_\alpha)$ contains $(I_F,b_\alpha|_F)$ where $I_F = I \cap F$.
\end{proposition}
\begin{proof}
Since $I$ is finitely generated over $\I_K$, it is then finitely generated over $\I_F$ as $\I_K$ is a finitely generated $\I_F$-module. Since $\I_F$ is a Noetherian ring, $I$ is Noetherian. Consider the set $I_F = I \cap F$ which is an $\I_F$-submodule of $F$ and an $\I_F$-submodule of $I$. Since $I$ is Noetherian, then $I_F$ is finitely generated, and is thus  a fractional ideal of $F.$  By definition of $I_F$, $(I_F,b_\alpha|_F)\subseteq (I,b_\alpha)$ is a $k$-subideal lattice of $(I,b_\alpha)$.

\end{proof}

We call the fractional ideal $I_F$ in Proposition \ref{Lemma: Existence of ideal below} the ideal below $I$ in $F.$ We note that if $I$ is an integral ideal in $\I_K$, then $I_F = I \cap F = I\cap \I_F$.
\begin{proposition}\label{Prop: upperbound of minimal discriminant of sublattice of rank k}
 Let $F$ be a subfield of $K$ of degree $k< n$ and consider the ideal lattice $(I,b_\alpha)$. The ideal $I_F$ below $I$ in $F$ provides the following upper bound on $d_k(I,b_\alpha):$
     \[ d_k(I,b_\alpha)\leq \operatorname{disc}(I_F,b_\alpha|_F) =
     \mathrm{N}_{F/\Q}(\mathrm{Tr}_{K/F}(\alpha))\mathrm{N}_{F/\Q}(I_F)^2 |\Delta_F| .\] Moreover, $(I_F,b_\alpha|_F)$ gives the tightest upper bound for $d_k(I,b_\alpha)$ among all $k$-subideal lattices in $F$.
    
\end{proposition}
\begin{proof}
   Consider the $k$-subideal lattice $(I_F,b_\alpha|_F)$ and let $\set{u_1,\ldots,u_k}$ be a $\Z$-basis of $I_F$. A generator matrix of the sublattice $(I_F,b_\alpha|_F)$ of $(I,b_\alpha)$ is given by 
   \[ M_F = \begin{pmatrix}
       \longleftarrow & \sigma_\alpha(u_1)& \longrightarrow \\
      &\vdots& \\
      \longleftarrow & \sigma_\alpha(u_k)& \longrightarrow \\
   \end{pmatrix}\] which is a matrix of size $k\times n$ and the corresponding Gram matrix is of the form $G_F = \left(\mathrm{Tr}_{K/\Q}(\alpha u_i\overline{u_j})\right)$ where $1\leq i,j \leq k$. If $K$ is CM, we consider the generator matrix above using the canonical complex embeddings. Note that
\[ \mathrm{Tr}_{K/\Q}(\alpha u_i \overline{u_j})= \mathrm{Tr}_{F/\Q}(\mathrm{Tr}_{K/F}(\alpha u_i\overline{u_j}))=\mathrm{Tr}_{F/\Q}\left(\mathrm{Tr}_{K/F}(\alpha) u_i\overline{u_j}\right)\] since $u_i,u_j \in I_F \subset F.$ Let $\mathrm{Emb}_\Q(F,\C)= \set{\eta_1,\ldots,\eta_k}$ be the set of embeddings of $F$ over $\Q$. The Gram matrix $G_F$ can be decomposed as follows:
\begin{equation*}
    \begin{split}
         G_F &= \left( \mathrm{Tr}_{F/\Q}\left(\mathrm{Tr}_{K/F}(\alpha) u_i\overline{u_j}\right) \right)_{1\leq i,j \leq k}\\
         &= \begin{pmatrix}
             \eta_1(u_1)& \ldots & \eta_k(u_1)\\
             \vdots & \ldots & \vdots\\
             \eta_1(u_i)& \ldots &\eta_k(u_i)\\
             \vdots & \ldots & \vdots\\
              \eta_1(u_k)& \ldots & \eta_k(u_k)\\
         \end{pmatrix}
         D
         \begin{pmatrix}
             \eta_1(\overline{u_1}) & \ldots &  \eta_1(\overline{u_k})\\
             \vdots & \ldots & \vdots\\
             \eta_i(\overline{u_1}) & \ldots &  \eta_i(\overline{u_k})\\
             \vdots & \ldots & \vdots\\
             \eta_k(\overline{u_1}) & \ldots &  \eta_k(\overline{u_k})\\
         \end{pmatrix}
    \end{split} 
\end{equation*} 
where $D=\mbox{diag}(\eta_1(\mathrm{Tr}_{K/F}(\alpha)), \ldots, \eta_k(\mathrm{Tr}_{K/F}(\alpha)))$. This decomposition allows to compute the sublattice discriminant $\operatorname{disc}(I_F,b_\alpha|F)=\det(G_F) $:
\begin{equation*}
    \begin{split}
        \det(G_F) &= \mathrm{N}_{F/\Q}(\mathrm{Tr}_{K/F}(\alpha)) \det\left((\eta_i(u_j))_{1\leq i,j\leq k}(\eta_i(\overline{u_j}))^\top_{1\leq i,j\leq k})\right)\\
        &= \mathrm{N}_{F/\Q}(\mathrm{Tr}_{K/F}(\alpha))\det\left((\mathrm{Tr}_{F/\Q}(u_i\overline{u_j}))_{1\leq i,j \leq k}\right)\\
        &= \mathrm{N}_{F/\Q}(\mathrm{Tr}_{K/F}(\alpha)) \mathrm{N}_{F/\Q}(I_F)^2 |\Delta_F|.
    \end{split}
\end{equation*} 

Lastly, note that $I_F$ is the largest fractional ideal in $F$ contained in $I$. So if $J \neq I_F$ is another fractional ideal in $F$ contained in $I,$ then $J \subsetneq I_F$ which means $\operatorname{disc}(J,b_\alpha|_F) > \operatorname{disc}(I_F,b_\alpha|_F).$
\end{proof}

\begin{corollary}\label{cor:gammaideal}
    Let $(I,b_\alpha)$ be an ideal lattice in $K$. For each $k<n$, then
    \[ \gamma^{ideal}_{n,k}(I,b_\alpha)= \min_{\substack{F\subseteq K\\ [F:\Q]= k}}\frac{\mathrm{N}_{F/\Q}(\mathrm{Tr}_{K/F}(\alpha))\mathrm{N}_{F/\Q}(I_F)^2 |\Delta_F| }{(\mathrm{N}_{K/\Q}(\alpha)\mathrm{N}_{K/\Q}(I)^2|\Delta_K|)^\frac{k}{n}}.\]
\end{corollary}
\begin{proof}
    Since $K$ is a number field of degree $n<\infty$, then there are only finitely many subfields $F$ of degree $k<n$. For each subfield $F$, $I_F$ gives the smallest discriminant among the $k$-subideal lattices by Proposition \ref{Prop: upperbound of minimal discriminant of sublattice of rank k}. So comparing the discriminant of $(I_{F},{b_\alpha}|_{F})$ across every subfield $F$ of degree $k$ and recalling that  
    \begin{eqnarray*}
    \gamma^{ideal}_{n,k}(I,b_\alpha) & = & \frac{\min_{(I_F,b_\alpha|_F)\subseteq (I,b_\alpha)}\operatorname{disc}(I_F,b_\alpha|_F)}{\operatorname{disc}(I,b_\alpha)^{k/n}}\\
    &=&
    \min_{\substack{F\subseteq K\\ [F:\Q]= k}}\frac{\mathrm{disc}(I_F,{b_\alpha}|_F)}{(\mathrm{N}_{K/\Q}(\alpha)\mathrm{N}_{K/\Q}(I)^2|\Delta_K|)^\frac{k}{n}}
    \end{eqnarray*}
    gives the claimed result, using the expression for $\mathrm{disc}(I_F,{b_\alpha}|_F)$ from Proposition \ref{Prop: upperbound of minimal discriminant of sublattice of rank k} .
\end{proof}

We provide two examples using the cyclotomic field $K=\Q(\zeta_{p^r})$ of degree $n=\varphi(p^r)=(p-1)p^{(r-1)}>2$ where $p$ is a prime, $r$ is a positive integer and $\zeta_{p^r}$ is a primitive $p$-th root of unity. It is known that $\operatorname{Gal}(K_{p^r}/\Q)\cong (\Z/p^r\Z)^\times$ which is cyclic if and only if $p$ is an odd prime. 

\begin{example}\label{ex:Kpr}
   If $p$ is odd, then any proper divisor of $n$ is $k=up^{j}$ where $u$ divides $(p-1)$, $j\leq r-1$ such that $\operatorname{gcd}(u,p)=1.$ For each $k$, there is a unique subfield $F$ of degree $k$ with absolute discriminant given by \cite[Theorem 4.1]{da2002computing}
    \[ |\Delta_F|=p^{u((j+2)p^j-\frac{p^{j+1}-1}{p-1})-1}.\] Let $\alpha = 1$ and consider the ideal lattice $(\I_{K_{p^r}},b_1).$ Then $\mathrm{N}_{F/\Q}(\mathrm{Tr}_{K/F}(1))=(n/k)^k$ and 
    \[ \gamma^{ideal}_{n,k}(\I_K,b_1)=\gamma_{n,k}((\I_K,b_1),(\I_F,b_1|_F))=(\frac{n}{k})^k\frac{|\Delta_F|}{|\Delta_K|^\frac{k}{n}}\] by Corollary \ref{cor:gammaideal}. In particular, if $k = \frac{n}{2}$, then $F= K^+$, the maximal real subfield of the cyclotomic field $K$ and since the absolute discriminants of the respective number fields are 
    \[ |\Delta_K| = p^{p^{r-1}(pr-r-1)},~ |\Delta_{K^+_{p^r}}|=p^{\frac{p-1}{2}((r+1)p^{r-1}-\frac{p^r-1}{p-1})-1},
    \] 
    this gives 
    \[ \gamma^{ideal}_{n,\frac{n}{2}}(\I_K,b_1)= \gamma_{n,\frac{n}{2}}((\I_K,b_1),(\I_{K^+},b_1|_{K^+}))=\frac{2^\frac{n}{2}}{\sqrt{p}}.\]
\end{example}

\begin{example}
For the case $K=\Q(\zeta_{2^r})$ where $r>2$, there could be multiple subfields of the same degree. We claim that for any divisor $k = 2^{j-1}$ where $1< j <r$, the minimum discriminant of a subfield of degree $k$ is attained by a cyclotomic field. Hence, for each proper divisor $k= 2^{j-1}$ of $n$ where $F$ is the unique cyclotomic subfield of degree $k$ such that $1<j<r$, we have 
    \[\gamma_{n,k}^{ideal}(\I_{K},b_1)=\gamma_{n,1}((\I_{K},b_1),(\I_{F},b_1|_F))= (2^{r-j})^{2^{j-1}}\frac{2^{2^{j-1}(j-1)}}{(2^{2^{r-1}(r-1)})^{2^{j-r}}}=1.\]

We now prove the claim.
    Suppose that $k=2^{j-1}$ for $1<j<r $ and let $F\subset K$ be of degree $k.$ By \cite[Theorem 3.2]{lopes2003discriminant}, 
\[ |\Delta_F|= \begin{cases}
    2^{2^{(j-1)}(j-1)}& \text{ if } F \text{ is a cyclotomic field,}\\
    2^{j2^{(j-1)}-1} & \text{ otherwise.}
\end{cases} \] It is clear that $2^{j-1}>1$ as $j>1$, so
\[2^{2^{(j-1)}(j-1)}= 2^{j2^{(j-1)}-2^{(j-1)}}<2^{j2^{(j-1)}-1}.\] 
The subfield which corresponds to a cyclotomic field thus attains the minimum absolute discriminant. 
\end{example}

Beside restricting to fractional ideals in $F$, a natural question is whether there exists any free $\Z$-module $M$ of rank $k$ that is not a fractional ideal such that $I_F \subsetneq M\subset I$. The existence of such an $M$ would give $d_k(I,b_\alpha)\leq \operatorname{disc}(\sigma_\alpha(M))< \operatorname{disc}(I_F,b_\alpha|_F)$, which is a better approximation. The next proposition shows that no such an $M$ exists.

\begin{proposition}
    There is no free $\Z$-module of rank $k$ in $I$ that contains $I_F.$
\end{proposition}
\begin{proof}
    Suppose $M$ is a free $\Z$-submodule of $I$ of rank $k$ such that $I_F \subseteq M \subset I$, so $|M/I_F|< \infty$ and hence, $M/I_F$ is torsion. However $I/I_F$ is torsion free. Indeed, let $x \in I$ and suppose that there exists $n \in \N$ such that $n(x+I_F)= nx+I_F = I_F \iff nx \in I_F.$ Since $I_F \subset F$, then $x = \frac{nx}{n} \in F$ and so, $x \in I \cap F = I_F$ which implies that $x+I_F = I_F$. So the only torsion element in $I/I_F$ is $I_F$, proving that $I/I_F$ is torsion free. Since $M/I_F \subseteq I/I_F$, the only torsion subgroup of a torsion free group $I/I_F$ is the trivial group, so $|M/I_F|=0 \iff M= I_F.$
\end{proof}

\subsection{The Case of $\gamma^{ideal}_{n,k}(\I_K,b_1)$}

In the particular case where $\alpha=1$ and $I=\mathcal{O}_K$, then the formula from Corollary \ref{cor:gammaideal} simplifies to
\[
\left(\frac{n}{k}\right)^k\min_{\substack{F\subset K\\ [F:\Q]= k}}\frac{ |\Delta_F| }{|\Delta_K|^\frac{k}{n}},
\] 
a quantity we bound next.

\begin{proposition}\label{prop:nk}
For any $2\leq k< n$, we have 
\[
\sup_{\substack{K\\ [K:\Q]=n}} \gamma_{n,k}^{ideal}(\I_K,b_1)\leq \left(\frac{n}{k}\right)^k
\] 
and equality holds exactly when $\min_{\substack{F\subset K\\ [F:\Q]=k}}\mathrm{N}_{F/\Q}(\Delta_{K/F})^{\frac{-k}{n}}= 1$. In particular, for $k=1$, we have
\[
\sup_{\substack{K\\ [K:\Q]=n}} \gamma_{n,1}^{ideal}(\I_K,b_1)= n \min_{\substack{K\\ [K:\Q]= n}} |\Delta_K|^{-\frac{1}{n}}.
\]
\end{proposition}
\begin{proof}
Suppose now $k>1$. 
It is known that \cite[Corollary 3.2.10]{neukirch2013algebraic}
\[ |\Delta_K|= \mathrm{N}_{F/\Q}(\Delta_{K/F})|\Delta_F|^{\frac{n}{k}}\] where $\Delta_{K/F}$ is the relative discriminant of $K$ with respect to $F$ which is an integral ideal in $F.$ Since $\Delta_{K/F}$ is an integral ideal in $\I_F$, then $\mathrm{N}_{F/\Q}(\Delta_{K/F})$ is a positive integer, hence $\mathrm{N}_{F/\Q}(\Delta_{K/F})^{\frac{-k}{n}} \leq 1.$ This gives 
\[ \left(\frac{n}{k}\right)^k\min_{\substack{F\subset K\\ [F:\Q]=k}}\frac{|\Delta_F|}{|\Delta_K|^\frac{k}{n}}= \left(\frac{n}{k}\right)^k\min_{\substack{F\subset K\\ [F:\Q]=k}}\mathrm{N}_{F/\Q}(\Delta_{K/F})^{\frac{-k}{n}} \leq \left(\frac{n}{k}\right)^k,\] showing that $\sup_{\substack{K\\ [K:\Q]=n}} \gamma_{n,k}^{ideal}(\I_K,b_\alpha)\leq \left(\frac{n}{k}\right)^k$.

If $k=1$, then $F=\mathbb{Q}$ and the bound given above is not tight since there are no non-trivial unramified extension over $\Q$. So
\[ 
\gamma^{ideal}_{n,1}(\I_K,b_1)=n |\Delta_K|^{-\frac{1}{n}}
\] and
\[
\sup_{\substack{K\\ [K:\Q]=n}} \gamma_{n,1}^{ideal}(\I_K,b_1)= n \min_{\substack{K\\ [K:\Q]= n}} |\Delta_K|^{-\frac{1}{n}}.
\] 

\
\end{proof}

The value of $\sup_{\substack{K\\ [K:\Q]=n}} \gamma_{n,1}^{ideal}(\I_K,b_1)$ can be approximated by knowing the minimum discriminant of the number field $K$ of a fixed degree $n$, which is still an open problem in algebraic number theory \cite{odlyzko1990bounds}. Some data on the minimum absolute discriminant from  \cite[Table 1]{odlyzko1990bounds} and \cite{y1987petits} can be found, see Table \ref{tab:smalldisc}, for number fields of degree up to 8, which serves as a best case approximation of using $\gamma^{ideal}_{n,1}(\I_K,b_1)$. Note that the field attaining the minimum discriminant for degree 3 to 7 is neither totally real or CM, hence serving as a strict upper bound while the case of $n=2$ and $8$ are both equality.

\begin{table}
\[
\begin{array}{c| c| c}

 n & \min_K|\Delta_K| & \text{upper bound for }n(\min_K|\Delta_K|)^{-\frac{1}{n}} \\
\hline
2 & 3 \text{ (CM)} & \frac{2}{3^{\frac{1}{2}}}\approx 1.11547  \\
3 & 23 & \frac{3}{23^\frac{1}{3}}\approx 1.0549 \\
4 & 117 & \frac{4}{117^\frac{1}{4}}\approx 1.21622  \\
5 & 1609 & \frac{5}{1609^\frac{1}{5}}\approx 1.14198 \\
6 & 9747 &  \frac{6}{9747^\frac{1}{6}}\approx 1.29819 \\
7 & 184607 & \frac{7}{184607^\frac{1}{7}}\approx 1.23816 \\
8 & 1257728 \text{ (CM)} & \frac{8}{1257728^{\frac{1}{8}}} \approx 1.38242
\end{array}
\]
\caption{\label{tab:smalldisc}}
Computation of an upper bound of $\sup_{\substack{K\\ [K:\Q]=n}} \gamma_{n,1}^{ideal}(\I_K,b_1)$, shown in the 3rd column.

\end{table}

We discuss below the achievability of the upper bound $(n/k)^k$ provided in the previous proposition on  $\sup_{\substack{K\\ [K:\Q]=n}} \gamma_{n,k}^{ideal}(\I_K,b_1)$. 

\begin{proposition}\label{Prop: min attained iff all subfield of same degree unramified}
    Let $K$ be a number field of degree $n>2$ and $k\vert n$ with $k>1$. Then \[\min_{\substack{F\subset K\\ [F:\Q]=k}}\mathrm{N}_{F/\Q}(\Delta_{K/F})^\frac{-k}{n} =1\] if and only if $K/F$ is unramified for all subfields $F$ of degree $k.$
\end{proposition}
\begin{proof}
   Suppose that 
   there exists a subfield $\Tilde{F}\subset K$ such that 
   $N_{\tilde{F}/\Q}(\Delta_{K/\tilde{F}})^\frac{-k}{n}= 1$. Then
   $\Delta_{K/\Tilde{F}}=\I_{\Tilde{F}} \iff K/\Tilde{F}$ is unramified. 
   Let $F\neq \Tilde{F}$ be any other subfield of degree $k$. By the transitivity of discriminants, 
\[|\Delta_{F}|^{\frac{n}{k}} \leq N_{F/\Q}(\Delta_{K/{F}})|\Delta_{F}|^{\frac{n}{k}}= |\Delta_K|=|\Delta_{\Tilde{F}}|^\frac{n}{k} \iff  |\Delta_{F}|\leq|\Delta_{\Tilde{F}}|.\] So
\[ \min_{\substack{F\subset K\\ [F:\Q]=k}}\frac{|\Delta_F|}{|\Delta_K|^\frac{k}{n}} = N_{\Tilde{F}/\Q}(\Delta_{K/\Tilde{F}})^\frac{-k}{n}= \frac{|\Delta_{\Tilde{F}}|}{|\Delta_K|^\frac{k}{n}} = \max_{\substack{F\subset K\\ [F:\Q]=k}}\frac{|\Delta_F|}{|\Delta_K|^\frac{k}{n}}\] if and only if  
\[\min_{\substack{F\subset K\\ [F:\Q]=k}}|\Delta_F|=\max_{\substack{F\subset K\\ [F:\Q]=k}} |\Delta_F|.\]  This means that all subfields of $K$ of degree $k$ have the same discriminant 
and the equality is achieved if and only if $N_{F/\Q}(\Delta_{K/F})= 1$ for all subfield $F \iff K/F$ is unramified for all $F$.
\end{proof}

 The next example shows that is it possible for $K$ to have a single subextension $F$ of a given degree (of degree 2 in the example) such that $K/F$ is unramified by the means of Hilbert class fields \cite[p.399]{neukirch2013algebraic}. In other words, the bound $(n/k)^k$ is achievable for some parameters $k,n$ (at least for $n=10$ and $k=2$).

\begin{example}
  Let $F=\mathbb{Q}(\sqrt{401})$ and let $K=F[X]/(X^5-X^4-5X^3+4X^2+3X-1)$. The field $K$ is the Hilbert class field of $F$ \cite{cohen2000computing}, and therefore the extension $K/F$ is unramified. Moreover, $K$ is totally real and uniquely contains the quadratic subfield $F$. It follows from Proposition \ref{Prop: min attained iff all subfield of same degree unramified} that 
  \[\gamma^{ideal}_{10,2}(\I_K,b_1) = (\frac{10}{2})^2 = 25.\]
\end{example}

\subsection{The Case of $\gamma_{n,1}$}
\label{subsec:gamman1}
As mentioned in Proposition \ref{Prop: upperbound of minimal discriminant of sublattice of rank k}, $\gamma_{n,k}(I,b_\alpha)\leq \gamma^{ideal}_{n,k}(I,b_\alpha)$. In this subsection, we introduce a family of lattices where equality is attained. This family of lattices is also known to be well-rounded \cite{fukshansky2012well}.

Set $K_m=\Q(\zeta_m)$ for $\zeta_m$ a primitive $m$-th root of unity, which is 
of degree $n= \varphi(m)$ for $\varphi$ the Euler totient function. Its discriminant is
\[
\Delta_{K_m}=
(-1)^{n/2}m^{n}\prod_{p\mid m} p^{-\frac{n}{p-1}}.
\]

\begin{proposition}\label{Prop: minimum vector attained by a rational number}
    Let $(I,b_\alpha)$ be an ideal lattice. If $d_1(I,b_\alpha)$ can be attained by some $r \in \Q,$ then $\gamma_{n,1}(I,b_\alpha)= \gamma^{ideal}_{n,1}(I,b_\alpha).$
\end{proposition}
\begin{proof}
    Without loss of generality, let $r\in \Q$ be positive such that $d_1(I,b_\alpha)= b_\alpha(r,r).$ So $r \in I\cap \Q = I_\Q$. Since $\Z$ is a PID, then $I_\Q = x\Z$ for some $x \in \Q$ positive, so $r = xt$ for some $t \in \Z$. Then we see that 
    \[ d_1(I,b_\alpha)= b_\alpha(r,r)= r^2n=(xt)^2n.\] If $t > 1$, then $(xt)^2= x^2t^2> x^2$ and so
    \[ b_\alpha(x,x)= x^2n<(xt)^2n= d_1(I,b_\alpha) \] which contradicts the minimality of $d_1(I,b_\alpha)$. Hence, $x\Z = r\Z$ and 
    \[ \gamma^{ideal}_{n,1}(I,b_\alpha)= \frac{\mathrm{disc}(r\Z,b_1|_\Q)}{\mathrm{disc}(I,b_\alpha)^\frac{1}{n}}=\frac{r^2}{\mathrm{disc}(I,b_\alpha)^\frac{1}{n}}= \gamma_{n,1}(I,b_\alpha).\]
\end{proof}
\begin{proposition}\label{prop:k2m}
The ideal lattice $(\I_{K_m},b_1)$, $m\geq 1$, satisfies
\[
\gamma_{n,1}(\I_{K_m},b_1)=\gamma^{ideal}_{n,1}(\I_{K_m},b_1)=n \frac{\prod_{p\vert m}p^{\frac{1}{p-1}}}{m}.
\]
\end{proposition}
\begin{proof}
We know that $d_1(\I_{K_m},b_1)$ is attained by the roots of unity \cite{fukshansky2012well} in $\I_{K_m}$ and that $\I_{K_m}\cap \Q = \Z$ which is generated by $1$. By Proposition \ref{Prop: minimum vector attained by a rational number}, then $\gamma^{ideal}_{n,1}(\I_{K_m},b_1)= \gamma_{n,1}(\I_{K_m},b_1).$ Finally, 
 \[ 
\gamma^{ideal}_{n,1}(\I_{K_m},b_1)= \frac{n}{|\Delta_{K_m}|^{\frac{1}{n}}}= \frac{n}{
\left(
m^{n}\prod_{p\mid m} p^{-\frac{n}{p-1}}
\right)^{1/n}}
=n\frac{\prod_{p\vert m}p^{\frac{1}{p-1}}}{m}.
 \]   
\end{proof}

General bounds for the Hermite constant are known \cite{cryptoeprint:2022/677}, namely
\[
\frac{n}{2\pi e}<\gamma_{n,1}< \frac{n}{8.5}+2,\]
which holds for every positive integer $n$. Since this estimate is uniform in $n$, it is generally not expected to be sharp for any particular dimension. We show below that for even dimensions at most 24 except 14 and 22, we can attain a better lower bound than $\frac{n}{2\pi e}$ by exhibiting a family $L_n=(I_n,b_1)$ of lattices such that 
\[
\frac{n}{2\pi e}<
\gamma_{n,1} (L_n)  \leq  \gamma_{n,1}.
\]

%
The family we will be using is formed of ideal lattices coming from the ring of integers of cyclotomic fields where $\alpha =1$. Consider the function
 $f(m) = \prod_{p\vert m}p^{\frac{1}{p-1}}.$ From Proposition \ref{prop:k2m}, we will be searching for values of $m$ such that  
 \[ \frac{f(m)}{m}\geq \frac{1}{2\pi e} \iff 2f(m)\pi e \geq m \geq \prod_{p\vert m}p \iff \prod_{p\vert m}p^{\frac{p-2}{p-1}}\leq 2\pi e.\] 
 
Suppose $m$ is a positive integer such that the exponents appearing in its prime factorization are all 1. 
Recall that $\Q(\zeta_m)=\Q(\zeta_{2m})$ for $m$ odd. 

\begin{itemize}
\item A single prime: $m=p=19$ is the largest prime such that $p^{\frac{p-2}{p-1}}\leq 2\pi e,$ hence the largest possible $m$ from this case is $m=2(19).$ Concretely, $p\in\{3,5,7,11,13,17,19\}$ and 
$2 p\in\{2( 3),2( 5),2( 7),2( 11),2( 13),2( 17), 2(19)\}$.
\item 
A product of two primes, one of which is $3$: $m=3(11)$ is the largest integer that satisfies the above inequality, hence the largest possible $m$ is $m=2 ( 3)(11).$ Concretely, $3 p\in\{3( 5),3( 7),3(11)\}$ and\\ $2( 3)(p)\in\{2( 3)(5), 2(3)( 7),2( 3)(11)\}$. 
\item 
A product of two primes, one of which is $5$: $m=5(7)$ and $m=2(5)(7)$ are the only two possibilities.
\end{itemize}

For the product of two primes, one of which is $>5$, or for the product of at least three primes,  no such an $m$ exists.  So we are looking at values of $m$ containing prime factors in the three cases above, as shown in the first column of Table  \ref{tab:improv} and the values in the third column are the values of $m$ by allowing exponents of prime to be greater than 1 such that $\gamma_{n,1}(\I_{K_m},b_1)\geq \frac{n}{2\pi e}$.
 
 \begin{table}
\[
\begin{array}{c|c|c|c|c}
\text{Factors} & 2f(m)\pi e\approx & m &n \\
\hline
2 &  34.1589 & 4,8,16,32&2,4,8,16\\
3 & 29.5825 & 3,9,27&2,6,18\\
5 &  25.5398 & 5,25&4,20\\
7 &  23.6224 & 7&6\\
11 &  21.7077 & 11&10\\
13&  21.1496 & 13&12\\
17 &  20.3881 &17&16\\
19 & 20.1148 & 19&18\\
2,3&59.1650&6,12,18,24,36,48&2,4,6,8,12,16\\
2,5&51.0795&10,20,40,50&4,8,16,20\\
2,7&47.2448&14,28&6,12\\
2,11&43.4154&22&10\\
2,13&42.2993&26,52&12,24\\
2,17&40.7763&34&16\\
2,19&40.2296&38&18\\
3,5&44.2362&15&8\\
3,7&40.9152&21&12\\
3,11&37.5988&33&20\\
2,3,5&88.4723&30,60&8,16\\
2,3,7&81.8304&42&12\\
2,3,11&75.1976&66&20\\
5,7&35.3237&35&24\\
2,5,7&70.6475&70&24
\end{array}\]
\caption{\label{tab:improv}
Parameters outperforming the classical lower bound on $\gamma_{n,1}$.
}
\end{table}

For all even dimension $2\leq n\leq 24$ except 14 and 22, the ideal Rankin invariants $\gamma^{ideal}_{n,1}(\I_{K_m},b_1)$ outperforms the lower bound $\frac{n}{2\pi e}.$

In \cite{cryptoeprint:2022/677}, it is conjectured that for dimension $n= 9$ to $23$, $\gamma_{n,1}=\frac{4}{\alpha_{4,n}^\frac{1}{n}}$ where $\alpha_{4,n}$ is the minimum discriminant of any $n$ dimension norm 4 integral lattice. Table \ref{tab:comparing computed,actual and conjectured value} compares for each dimension $n$, the largest value of $\frac{nf(m)}{m}$ attained corresponding to a value of $m$ with $\gamma_{n,1}$ and $\frac{4}{\alpha_{4,n}^\frac{1}{n}}.$
\begin{table}
\[
\begin{array}{c|c|c|c|c|c}
n&m&\frac{nf(m)}{m}\approx&\gamma_{n,1}\approx&\frac{4}{\alpha_{4,n}^\frac{1}{n}} & \text{upper bound for }n(\min_K|\Delta_K|)^{-\frac{1}{n}}  \\
\hline
2&3&1.1547&1.1547&1.1547 & \approx 1.1547\\
4&5&1.1963&1.4142&1.4142 & \approx 1.21622\\
6&7&1.1855&1.6654&1.6654 & \approx 1.29819 \\
8&15&1.3813&2&2&\approx 1.38242\\
10&11&1.1554&&2.0584\\
12&28&1.1855&&2.3094\\
16&48&1.1547&&2.8284\\
18&27&1.1547&&2.9868\\
20&33&1.3342&&3.249\\
24&35&1.4182&&4
\end{array}\]
\caption{\label{tab:comparing computed,actual and conjectured value}
Comparison between $\frac{nf(m)}{m}, \gamma_{n,1}$, $\frac{4}{\alpha_{4,n}^\frac{1}{n}}$ and the upper bound given by the smallest discriminant. 
}
\end{table}

\subsection{The Role of Primes}

For a number field $K$, let $J_K$ be the set of fractional ideals of $K$ and for each prime ideal $\mathfrak{P}$ in $\I_K$, we define the $\mathfrak{P}$-adic valuation on $J_K$ to be the map 
\begin{equation*}
    \begin{split}
        v_{\mathfrak{P}}: J_K &\rightarrow \Z\\
        I &\mapsto v_{\mathfrak{P}}(I)
    \end{split}
\end{equation*} where $v_{\mathfrak{P}}(I)$ is the largest integer such that $\mathfrak{P}^{v_{\mathfrak{P}}(I)}$ is in the prime factorization of $I.$ For any fractional ideal $I\in J_K$, we then have 
\[ I = \prod_{\mathfrak{P}}\mathfrak{P}^{v_{\mathfrak{P}}(I)}\] where almost all but finitely $v_{\mathfrak{P}}(I)\neq 0$. 

Let $K/F$ be a field extension and $\mathfrak{P}$ be a prime ideal in $\I_K$ with $\mathfrak{p}= \mathfrak{P}\cap \I_F$. Denote respectively the ramification index and the inertial degree of $\mathfrak{P}$ over $\mathfrak{p}$ by $e_{\mathfrak{P}|\mathfrak{p}}= v_{\mathfrak{P}}(\mathfrak{p}\I_K)$
and $f_{\mathfrak{P}|\mathfrak{p}} = [\I_K/\mathfrak{P}:\I_F/\mathfrak{p}]$. In the case where $K/F$ is Galois, then all the primes $\mathfrak{P}|\mathfrak{p}$ have the same ramification index and inertial degree and we denote it by $e_\mathfrak{p}$ and $f_{\mathfrak{p}}$ respectively. 

\begin{proposition}\label{Prop: factorization of valuation}
    For any fractional ideal $I \in J_F$, we have 
    \[ v_{\mathfrak{P}}(I\I_K)= e_{\mathfrak{P}|\mathfrak{p}}v_{\mathfrak{p}}(I).\]
\end{proposition}
\begin{proof}
    Let $I = \prod_{\mathfrak{p}}\mathfrak{p}^{v_{\mathfrak{p}}(I)}$ be a fractional ideal in $F.$ Then 
    \[ I\I_K = (\prod_{\mathfrak{p}}\mathfrak{p}^{v_\mathfrak{p}(I)})\I_K= \prod_{\mathfrak{p}}(\mathfrak{p}\I_K)^{v_\mathfrak{p}(I)}= \prod_{\mathfrak{p}}\prod_{\mathfrak{P}|\mathfrak{p}}(\mathfrak{P}^{e_{\mathfrak{P}|\mathfrak{p}}})^{v_\mathfrak{p}(I)}
    ,\] hence, $v_\mathfrak{P}(I)= e_{\mathfrak{P}|\mathfrak{p}}v_{\mathfrak{p}}(I).$
\end{proof}

Let $I\in J_K$, since $I =\prod_{\mathfrak{P}}\mathfrak{P}^{v_{\mathfrak{P}}(I)} = \prod_{\mathfrak{p}} \prod_{\mathfrak{P}| \mathfrak{p}}\mathfrak{P}^{v_{\mathfrak{P}}(I)}$, we label the distinct primes below from $\mathfrak{p}_1$ to $\mathfrak{p}_t$, and for each $i$, we similarly label the distinct primes $\mathfrak{P}_{i,j}|\mathfrak{p}_i$ for $j=1,\ldots, r_i$ to get 
\[
I= \prod_{i=1}^t\prod_{j=1}^{r_i}\mathfrak{P}_{i,j}^{v_{\mathfrak{P}_{i,j}(I)}},~v_{\mathfrak{P}_{i,j}(I)}\in \Z.
\]

\begin{proposition}\label{Prop: general factorization of ideal below}
    We have
    \[ I_F  = \prod_{i=1}^t \mathfrak{p}_i^{m_{i}}\] where $m_{i}=\max_{j=1,\ldots,r_i} \left\lceil \frac{v_{\mathfrak{P}_{i,j}(I)}}{e_{\mathfrak{P}_{i,j}|\mathfrak{p}_i}}\right\rceil$ for each $i=1,\ldots, t$.
\end{proposition}
\begin{proof}

Write $I_F = \prod_{\mathfrak{p}}\mathfrak{p}^{v_\mathfrak{p}(I_F)}$. We note that $I_F \subseteq I \implies I_F \I_K \subseteq I$. By Proposition \ref{Prop: factorization of valuation}, for each $i=1,\ldots, t$ and $j=1,\ldots,r_i$, then 
  \[ v_{\mathfrak{P}_{i,j}}(I)\leq v_{\mathfrak{P}_{i,j}}(I_F\I_K) =e_{\mathfrak{P}_{i,j}|\mathfrak{p}_i}v_{\mathfrak{p}_i}(I_F) \iff v_{\mathfrak{p}_i}(I_F) \geq  \frac{v_{\mathfrak{P}_{i,j}(I)}}{e_{\mathfrak{P}_{i,j}|\mathfrak{p}_i}}. \] Since for each $i$, it holds for all $j$, we have that  
  $ v_{\mathfrak{p}_i}(I_F) \geq \max_{j=1,\ldots,r_i} \frac{v_{\mathfrak{P}_{i,j}(I)}}{e_{\mathfrak{P}_{i,j}|\mathfrak{p}_i}}. $ As $v_{\mathfrak{p}_i}(I_F)\in \Z$, then 
  \[  v_{\mathfrak{p}_i}(I_F) \geq \max_{j=1,\ldots,r_i} \left\lceil \frac{v_{\mathfrak{P}_{i,j}(I)}}{e_{\mathfrak{P}_{i,j}|\mathfrak{p}_i}}\right\rceil=m_i.\]
   This implies that $\mathfrak{p}_i^{{v_{\mathfrak{p}}}_i(I_F)} \subseteq \mathfrak{p}^{m_i}_i$ for all $i=1,\ldots, t.$ 
    Hence $I_F=\prod_{\mathfrak{p}}\mathfrak{p}^{v_\mathfrak{p}(I_F)}\subseteq \prod_{i=1}^t\mathfrak{p}^{m_i}_i$.
    
    Conversely, let $J = \prod_{i=1}^t\mathfrak{p}_i^{m_i}\subseteq F$, then by Proposition \ref{Prop: factorization of valuation}, for any $i,j,$
    \begin{equation*}
        \begin{split}
            v_{\mathfrak{P}_{i,j}}(J\I_K)&= e_{\mathfrak{P}_{i,j}|\mathfrak{p}_i}v_{\mathfrak{p}_i}(J)\\ &= e_{\mathfrak{P}_{i,j}|\mathfrak{p}_i} 
            m_i\\
            &\geq e_{\mathfrak{P}_{i,j}|\mathfrak{p}_i} 
            \left\lceil \frac{v_{\mathfrak{P}_{i,j}(I)}}{e_{\mathfrak{P}_{i,j}|\mathfrak{p}_{i}}}\right\rceil\\
            &\geq e_{\mathfrak{P}_{i,j}|\mathfrak{p}_i} 
             \frac{v_{\mathfrak{P}_{i,j}(I)}}{e_{\mathfrak{P}_{i,j}|\mathfrak{p}_{i}}}\\
            &= v_{\mathfrak{P}_{i,j}}(I).
        \end{split}
    \end{equation*} Hence $J\I_K \subseteq I $ and so,  $J\subseteq  J\I_K \cap \I_F \subseteq I \cap F\ $ that is $ J \subseteq I_F.$
\end{proof}
   
\begin{theorem}\label{thm:ideal gamma}
    Let $(I,b_\alpha)$ be an ideal lattice for $I=\prod_{i=1}^{t}\prod_{j=1}^{r_i}\mathfrak{P}_{i,j}^{v_{\mathfrak{P}_{i,j}(I)}} \in J_K$. The ideal Rankin invariant of order $k$ is 
\[
  \gamma^{ideal}_{n,k}(I,b_\alpha)= \min_{\substack{F\subseteq K\\ [F:\Q]= k}}\frac{\mathrm{N}_{F/\Q}(\mathrm{Tr}_{K/F}(\alpha)) }{(\mathrm{N}_{K/\Q}(\alpha)N_{F/\mathbb{Q}}(\Delta_{K/F})^\frac{k}{n}}\prod_{i=1}^t \frac{p_i^{2m_{i}f_{\mathfrak{p}_i|p_i}}}{p_i^{\frac{2k}{n}\sum_{j=1}^{r_i}v_{\mathfrak{P}_{i,j}(I)}f_{\mathfrak{P}_{i,j}|p_i}}} 
\]
for $p_i$ the prime below $\mathfrak{p}_i$ and  $m_{i}= \max_{j=1,\ldots,r_i}\left\lceil \frac{v_{\mathfrak{P}_{i,j}(I)}}{e_{\mathfrak{P}_{i,j}|\mathfrak{p}_i}}\right\rceil$ for each $i=1,\ldots,t$ and $j=1,\ldots, r_i.$             
\end{theorem}
\begin{proof}
Let $I = \prod_{i=1}^t\prod_{j=1}^{r_i}\mathfrak{P}^{v_{\mathfrak{P}_{i,j}(I)}}_{i,j}$ be the factorization of $I \in J_K$ where for each $i$, $\mathfrak{p}_i$ is the prime below $\mathfrak{P}_{i,j}$ for all $j=1,\ldots,r_i$.  For each $i=1,\ldots,t$, let $p_i\in \Z$ be the prime below $\mathfrak{p}_i$ i.e., $\mathfrak{p}_i\cap \Z = p_i$, then \[ \mathfrak{P}_{i,j}\cap \Z = \mathfrak{P}_{i,j}\cap \I_F\cap \Z = \mathfrak{p}_i \cap \Z = p_i.\]  
    By Proposition \ref{Prop: general factorization of ideal below}, 
     \[ \mathrm{N}_{F/\Q}(I_F)= \mathrm{N}_{F/\Q}( \prod_{i=1}^{t}\mathfrak{p}^{m_{i}}_i) =
     \prod_{i=1}^t\mathrm{N}_{F/\Q}(\mathfrak{p}_i)^{m_{i}}= \prod_{i=1}^tp_i^{m_{i}f_{\mathfrak{p}_i|p_i}}\] 
     where $m_{i} = \max_{j=1,\ldots,r_i}\left\lceil \frac{v_{\mathfrak{P}_{i,j}(I)}}{e_{\mathfrak{P}_{i,j}|\mathfrak{p}_i}}\right\rceil$ for each $i=1,\ldots, t.$ Furthermore
    \[ \mathrm{N}_{K/\Q}(I)= \prod_{i=1}^t\prod_{j=1}^{r_i}\mathrm{N}_{K/\Q}(\mathfrak{P}_{i,j})^{v_{\mathfrak{P}_{i,j}(I)}}= \prod_{i=1}^tp_i^{\sum_{j=1}^{r_i}v_{\mathfrak{P}_{i,j}(I)}f_{\mathfrak{P}_{i,j}|p_i}}.\]
Lastly, use $\frac{|\Delta_F|}{|\Delta_K|^\frac{k}{n}}= \mathrm{N}_{F/\Q}(\Delta_{K/F})^{\frac{-k}{n}}$ with the formula for $\gamma^{ideal}_{n,k}(I,b_\alpha)$ from Corollary \ref{cor:gammaideal}.

\end{proof}

\begin{example}\label{Ex: rank 10}
    Consider $K=\Q(\zeta_{25}+\zeta_{25}^{-1})$ which is a cyclic extension over $\Q$ of degree $10$ and $F= \Q[x]/(X^2+5X-25)$ which is the unique subfield of degree $2$. The different ideal of $K$ is $\mathfrak{D}_K= \mathfrak{P}^{17}$ where $\mathfrak{P}$ is the unique prime above $5$ which is totally ramified in $K$. In fact, $5$ is totally ramified in both $F$ and $K$. Set $I = \mathfrak{P}^{-8}$ and consider the ideal lattice $(I,b_\alpha)$ where $\alpha =1$ so that $\mathrm{disc}(I,b_\alpha)=5$. A Gram matrix is given by \[\begin{psmallmatrix}10 & 5 & 0 & 2 & 8 & 0 & 0 & 8 & 2 & 0 \\
5 & 3 & 0 & 1 & 4 & 0 & 0 & 4 & 1 & 0 \\
0 & 0 & 20 & 16 & 4 & 5 & 5 & 4 & 16 & 15 \\
2 & 1 & 16 & 14 & 5 & 4 & 4 & 5 & 13 & 12 \\
8 & 4 & 4 & 5 & 8 & 1 & 1 & 7 & 5 & 3 \\
0 & 0 & 5 & 4 & 1 & 2 & 1 & 1 & 4 & 4 \\
0 & 0 & 5 & 4 & 1 & 1 & 2 & 1 & 4 & 4 \\
8 & 4 & 4 & 5 & 7 & 1 & 1 & 8 & 5 & 3 \\
2 & 1 & 16 & 13 & 5 & 4 & 4 & 5 & 14 & 12 \\
0 & 0 & 15 & 12 & 3 & 4 & 4 & 3 & 12 & 12\end{psmallmatrix}\] of determinant $5$. Let $\mathfrak{p}$ be the unique prime below $\mathfrak{P}$ in $F$ with ramification index $e_\mathfrak{p}=5$. By Proposition \ref{Prop: general factorization of ideal below}, $I_F = \mathfrak{p}^{\left\lceil \frac{-8}{5}\right\rceil}= \mathfrak{p}^{-1}$ and using Theorem \ref{thm:ideal gamma}, we have that 
    \begin{equation*}
        \begin{split}
            \gamma_{10,2}^{ideal}(I,b_1)&=\frac{N_{F/\Q}(Tr_{K/F}(1))}{N_{K/\Q}(1)^\frac{2}{10}}N_{F/\Q}(15625\I_F)^{\frac{-2}{10}}\frac{5^{2(-1)(1)}}{(5^{\frac{2(2)}{10}}(-8)(1))}\\
            &= 5^2(5^{12})^{\frac{-2}{10}}(5)^{-2+\frac{16}{5}}\\
            &=5^\frac{4}{5}.
        \end{split}
    \end{equation*} If $\gamma_{10,2}(I,b_\alpha)= \gamma_{10,2}^{ideal}(I,b_\alpha)$, this implies that $d_2(I,b_\alpha)= 5$ which is not the case since from the Gram matrix of $(I,b_\alpha)$, we can find a Gram matrix of a sublattice that is $\begin{pmatrix}
        2&1\\1&2
    \end{pmatrix}$ of determinant $3$. 
    
    In fact, we can use the structure of the number field $K$ and the fact that $(I,b_1)$ is an integral lattice to prove that $d_2(I,b_1) =3.$ Let $G=\langle\tau \rangle$ be the Galois group of $K/\Q$. It suffices to show that $d_1(I,b_1)=2 $ and using Proposition \ref{Prop: general bounds on d_2(L)}, it gives $\gamma_{10,2}(I,b_1)\geq \frac{3}{4}d_1(I,b_1)^2 = 3$ which is attained by the sublattice mentioned above. Since $(I,b_1)$ is an integral lattice, then $d_1(I,b_1)\geq 1$. Suppose that there exists $x \in I$ such that $b_1(x,x)= 1.$ Since $K$ is Galois, for any $0\leq j\leq 9$, $\tau^j(\mathfrak{P})= \mathfrak{P}$ as $\mathfrak{P}$ is the unique prime above $5$, so $\tau^j(I)= I$ and $\tau^j(x) \in I.$ Note that for any $j$, 
    \[b_1(\tau^j(x),\tau^j(x))=\mathrm{Tr}_{K/\Q}(\tau^j(x^2))= \sum_{i=0}^{9}\tau^{i+j}(x^2)= \mathrm{Tr}_{K/\Q}(x^2)= b_1(x,x)=1.\] If for any distinct $i,j \in \{0,\ldots,9\}$, $b_1(\tau^i(x),\tau^j(x))= 0$, then the image of $\{x,\tau(x),\ldots,\tau^9(x)\}\subset I$ under $\sigma_\alpha$ spans a sublattice $L\cong \Z^{10}$. However,
    \[ \mathrm{disc}(L)= |(I,b_1)/L|^2\mathrm{disc}(I,b_1) \iff |(I,b_1)/L|^2 = \frac{1}{5} \not\in \Z\] which is a contradiction. Hence, there exists $i \neq j$ such that $b_\alpha(\tau^i(x),\tau^j(x))= 1$ and without loss of generality, let $i =0$, so $\tau^i(x)=x$. Then $\tau^j(x)= \pm x.$ If $\tau^j(x)=x$, then $x^2$ is fixed by the subgroup of $\langle\tau^j\rangle \subset G$, hence $x^2 \in K^{\langle \tau^j\rangle } \subsetneq K$ as $j\neq 0.$ If $\tau^j(x)= -x,$ then $ \tau^j(x^2)= \tau^j(x)^2 =x^2$ which implies that $x^2 \in K^{\langle\tau^j\rangle} \subsetneq K.$ Since $x \in I= \mathfrak{P}^{-8}\subseteq \mathfrak{P}^{-17}= \mathfrak{D}_K^{-1}\implies x^2\in \mathfrak{P}^{-16}\subseteq \mathfrak{P}^{-17}=\mathfrak{D}_K^{-1}$, and in either cases, $x^2\in K^{\langle\tau^j\rangle} \cap \mathfrak{P}^{-16}\subseteq K^{\langle \tau^j\rangle}\cap \mathfrak{D}^{-1}_K= \mathfrak{D}_{K^{\langle\tau^j\rangle}}^{-1}.$ This implies that $\mathrm{Tr}_{K^{\langle\tau^j\rangle}/\Q}(x^2)$ is a positive integer and 
    \[ b_\alpha(x,x)= Tr_{K/\Q}(x^2)= [K:K^{\langle \tau^j\rangle}]\mathrm{Tr}_{K^{\langle\tau^j\rangle}/\Q}(x^2) \geq 2,\] which contradicts the fact that we assumed $b_1(x,x)=1$. Therefore, we see that $b_1(x,x)\geq 2$ which is attained by one of the $\Z$-basis elements of the sublattice of discriminant 3. This proves that $d_1(I,b_1)=2$, hence this implies that $d_2(I,b_1)=3.$ Finally, 
    \[ \gamma_{10,2}\geq \gamma_{10,2}(I,b_1)= \frac{3}{5^\frac{2}{10}}\approx  2.1743.\]

    Next, for the case of $k=5$, we can also consider the unique subfield $F=\Q(y)=\Q[X]/(X^5 - 10X^3 + 5X^2 + 10X + 1)$ of degree $5$ and let $K = F(z)$. Similarly, we get that 
    \[\gamma^{ideal}_{10,5}(I,b_\alpha)= \frac{32}{5^\frac{1}{2}}\approx 14.3108.\] A $\Z$-basis for $I_F$ is $\{1,\frac{1}{5}(y+1),\frac{1}{35}(y^4+3y^3-8y^2-19y+16), \frac{1}{35}(-3y^4-2y^3+24y^2-6y+1), \frac{1}{35}(-2y^4+y^3+23y^2-18y-18)\}$. This is in fact not $\gamma_{10,5}(I,b_\alpha)$ since 1 with the element $(\frac{1}{7}(-\frac{1}{5}y^4+\frac{4}{5}y^3+3y^2-6y+1)z+\frac{1}{7}(\frac{1}{5}y^4+\frac{3}{5}y^3-\frac{1}{5}y^2-1-\frac{12}{5}))\in I$, we get the follow sublattice of rank 5 with Gram matrix 
    \[\begin{psmallmatrix}
        2&-1&0&0&1\\
        -1&2&0&0&0\\
        0&0&2&0&0\\
        0&0&0&2&0\\
        1&0&0&0&2
    \end{psmallmatrix} \] of determinant 16. So, 
    \[ \gamma_{10,5}(I,b_1)\leq \frac{16}{5^\frac{1}{2}}\approx 7.1554.\] However, to determine $d_5(I,b_1)$ is still a difficult task.
\end{example}

%
%

\section{On the Looseness of Subideal Lattices}
\label{sec:loose}

Given an ideal lattice $(I,b_\alpha)$, we expect $\gamma^{ideal}_{n,k}(I,b_\alpha)$  to be loose with respect to $\gamma_{n,k}(I,b_\alpha)$ since the minimum for the former is only computed over a subset of sublattices given by $I_F=I \cap F$ for $F$ ranging over subfields of $K$ of degree $k$ which is demonstrated in Example \ref{Ex: rank 10}. We will prove necessary conditions for a $k$-subideal lattice to attain optimality  We also mention particular examples where $k$-subideal lattices are attaining optimality.

We write $v_p$ for the valuation in $p$ of an integer, and start with the case where $k=1.$

\begin{proposition}\label{Prop: Necessary condition for subideal to attain optimality rank 1}
    Let $(I,b_\alpha)$ be an ideal lattice and $d_1(I,b_\alpha)= \prod_{i=1}^t p_i^{c_i} \in \Z$ such that $p_i$ are distinct primes and $c_i$ are nonnegative integers. If $v_{p_i}(\mathrm{Tr}_{K/\Q}(\alpha))\not\equiv c_i \pmod 2$ for some $i\in\{1,\ldots,t\}$ or $\frac{\mathrm{Tr}_{K/\Q}(\alpha)}{\prod_{i=1}^t p_i^{v_{p_i}(\mathrm{Tr}_{K/\Q}(\alpha))}}$ is not a square, then $\gamma_{n,1}(I,b_\alpha)\neq \gamma^{ideal}_{n,1}(I,b_\alpha)$.
\end{proposition}
\begin{proof}
    Suppose that $\gamma_{n,1}(I,b_\alpha)= \gamma^{ideal}_{n,1}(I,b_\alpha)$, so $d_1(I,b_\alpha)= \operatorname{disc}(I_\Q,b_\alpha) = \mathrm{Tr}_{K/\Q}(\alpha)r^2$ where $I_\Q = (r)$ for some $r \in \Q^\times$. Factorizing both $\mathrm{Tr}_{K/\Q}(\alpha)$ and $r$, for all $i=1,\ldots,t$, there exist $b_1,b_2 \in \Q$ such that $v_{p_i}(b_1)=v_{p_i}(b_2)=0$ and $Tr_{K/\Q}(\alpha)=(\prod_{i=1}^t p_i^{v_{p_i}(\mathrm{Tr}_{K/\Q}(\alpha))})b_1$ and $r= (\prod_{i=1}^tp_i^{v_{p_i}(r)})b_2$ and so, 
    \[ \prod_{i=1}^t p_i^{c_i} = (\prod_{i=1}^t p_i^{v_{p_i}(\mathrm{Tr}_{K/\Q}(\alpha))+2v_{p_i}(r)})b_1b_2^2 \iff c_i = v_{p_i}(\mathrm{Tr}_{K/\Q}(\alpha))+2v_{p_i}(r),~ b_1=b_2^{-2}.\] 
    This means that for all $i,$ $v_{p_i}(\mathrm{Tr}_{K/\Q}(\alpha))\equiv c_i \pmod 2$ and $b_1= \frac{\mathrm{Tr}_{K/\Q}(\alpha)}{\prod_{i=1}^tp_i^{v_{p_i}(\mathrm{Tr}_{K/\Q}(\alpha))}}$ is a square.
\end{proof}

Next, we look into the case $k>1$ where we assume that $F$ is a Galois extension over $\mathbb{Q}$.  We state an extension of \cite[Lemma 4.2]{jorge2013rotated} which will be useful.

\begin{lemma}\label{Lemma: Norm of ideal in galois extension}
   Let $F$ be a Galois extension over $\Q$ (not necessarily CM or totally real). Let $f_{p_i}$ be the inertial degree for any prime above the distinct primes $p_i\in \Z$ where $i=1,\ldots, t$. 
   \begin{enumerate}
       \item For any integral ideal $J\subseteq \I_K$, then 
   \[\mathrm{N}_{F/\Q}(J) = (\prod_{i=1}^t(p^{f_{p_i}}_i)^{a_i})b\] where $a_i \geq 0$ and positive integer $b$ such that $\operatorname{gcd}(b,p_i)=1$ for all $i.$
   \item Furthermore
   \[
   |\Delta_F|= (\prod_{i=1}^t(p^{f_{p_i}}_i)^{a_i})b
   \]
   where $a_i \geq 0$ and positive integer $b$ such that $\operatorname{gcd}(b,p_i)=1$ for all $i.$
   \end{enumerate}

\end{lemma}
\begin{proof}
  \begin{enumerate}
       \item Fix $p_1,\ldots, p_t$ to be distinct primes and let $J = \prod_{j=1}^r \mathfrak{q}_j^{v_{\mathfrak{q}_j}(J)}$ be the unique factorization of $J$ in $\I_F$. If for each $i=1,\ldots,t$, $p_i$ is not below $\mathfrak{q}_j$ for any $j$, then $\mathrm{N}_{F/\Q}(\mathfrak{q}_j^{v_{\mathfrak{q}_j}(J)})= q_j^{f_{\mathfrak{q}_j}v_{\mathfrak{q}_j}(J)}$ where $q_j$ are the primes below $\mathfrak{q}_j$ for each $j$ and $\mathrm{gcd}(q_j,p_i)=1$ for all $i,j.$ Set 
       \[b =\prod_{\substack{\mathfrak{q}_j \cap \Z \neq p_i\\ \forall i=1,\ldots, t}}\mathrm{N}_{K/\Q}(\mathfrak{q}_j^{v_{\mathfrak{q}_j}(J)})=\prod_{j=1}^r q_j^{f_{\mathfrak{q}_j}v_{\mathfrak{q}_j}(J)},\] then $\mathrm{gcd}(b,p_i)=1$ for all $i=1,\ldots,t$. Next, suppose there exists $i$ such that $p_i$ is below $\mathfrak{q}_j$ for some $j$, so $\mathrm{N}_{F/\Q}(\mathfrak{q}_j^{v_{\mathfrak{q}_j}(J)})= p_i^{f_{p_i}v_{\mathfrak{q}_j}(J)}$. For each $i$, 
       \[ \prod_{\mathfrak{q}_j\cap \Z = p_i} \mathrm{N}_{F/\Q}(\mathfrak{q}_j^{v_{\mathfrak{q}_j}(J)})= p_i^{f_{p_i}(\sum_{\mathfrak{q}_j\cap \Z = p_i}v_{\mathfrak{q}_j}(J))}\] and if we let $a_i = \sum_{\mathfrak{q}_j\cap \Z = p_i}v_{\mathfrak{q}_j}(J)$, then 
       \begin{equation*}
           \begin{split}
               \mathrm{N}_{F/\Q}(J)&=\prod_{j=1}^r\mathrm{N}_{F/\Q}(\mathfrak{q}_j^{v_{\mathfrak{q}_j}(J)})\\&= (\prod_{i=1}^t  \prod_{\mathfrak{q}_j\cap \Z = p_i} \mathrm{N}_{F/\Q}(\mathfrak{q}_j^{v_{\mathfrak{q}_j}(J)}))\prod_{\substack{\mathfrak{q}_j\cap \Z \neq p_i\\ \forall i=1,\ldots,t}} \mathrm{N}_{F/\Q}(\mathfrak{q}_j^{v_{\mathfrak{q}_j}(J)})\\ &= (\prod_{i=1}^t(p_i^{f_i})^{a_i})b.
           \end{split}
       \end{equation*}
\item
     This follows immediately from 1. since $|\Delta_F|$ is the norm of the different ideal of $F$. 
\end{enumerate}
\end{proof}

\begin{proposition}\label{Prop: Necessary condition for subideal lattice to attain optimality}
    Suppose that $K$ contains a subfield $F$ that is Galois over $\Q$ of degree $k>1$, $\alpha \in K$ is totally positive and let $(I,b_\alpha)$ be an ideal lattice where $I$ is a fractional ideal of $K$ and $d_k(I,b_\alpha)=\prod_{i=1}^tp_i^{c_i}$ such that $p_i$ are distinct primes and  $c_i$ are non-negative integers. Set $f_{p_i}$ to be the inertial degree of $p_i$ in $F$. If  $\gamma^{ideal}_{n,k}(I,b_\alpha)= \gamma_{n,k}(I,b_\alpha)$ is attained at $F$, then $f_{p_i} \vert c_i$ for all $i$.
    
    Moreover, suppose that $\mathrm{N}_{F/\Q}(\mathrm{Tr}_{K/F}(\alpha))$ is a square and $c_i$ are all positive odd integers. 
    If $\gamma^{ideal}_{n,k}(I,b_\alpha)= \gamma_{n,k}(I,b_\alpha)$ is attained at $F$, then $|\Delta_F|= b\prod_{i=1}^tp_i^{a_if_{p_i}}$ where $a_if_{p_i}$ is odd and $b$ is a square. In particular, $p_i$ must ramify in $F$ and $|\Delta_F|$ must be divisible by a square coprime to $p_i$ for all $i$.
\end{proposition}
\begin{proof}
 Suppose that $\gamma^{ideal}_{n,k}(I,b_\alpha)=\gamma_{n,k}(I,b_\alpha)$ is attained at $F$. Then
    \[ \prod_{i=1}^tp_i^{c_i}= \operatorname{disc}(I_F,b_\alpha|_F)= \mathrm{N}_{F/\Q}(\mathrm{Tr}_{K/F}(\alpha))\mathrm{N}_{F/\Q}(I_F)^2|\Delta_F|.\]  
    We apply Lemma \ref{Lemma: Norm of ideal in galois extension} on each of the three product terms.

    We start with $\mathrm{N}_{F/\Q}(I_F)$. 
    There exists $c_1\in \I_F$ such that $c_1 I_F \subseteq \I_F$ and in particular, since $c_1| \mathrm{N}_{F/\Q}(c_1) \in\Z$, then $\mathrm{N}_{F/\Q}(c_1)I_F = J \subseteq \I_F$, so we can further assume $c_1 \in \Z$ positive. Furthermore, \[\mathrm{N}_{F/\Q}(c_1) =c_1^k=|c_1|^k =\mathrm{N}_{F/\Q}(c_1\I_F),\] hence, 
    \[\mathrm{N}_{F/\Q}(I_F)= \frac{\mathrm{N}_{F/\Q}(J)}{\mathrm{N}_{F/\Q}(c_1\I_F)}=\prod_{i=1}^tp_i^{(a_{i1}-a_{i2})f_{p_i}}b_1b_2^{-1}\] where $a_{i1},a_{i2}\geq 0$ and positive integers $b_1,b_2$ such that $\operatorname{gcd}(b_1,p_i)=\operatorname{gcd}(b_2,p_i)=1$ for all $i.$ 
    
    Next, we consider $\mathrm{N}_{F/\Q}(\mathrm{Tr}_{K/F}(\alpha))>0.$ Let $c_2 \in \Z$ positive such that $c_2\mathrm{Tr}_{K/F}(\alpha)\in \I_F$, so there exist $a_{i3},a_{i4} \geq 0$ and positive integers $b_3,b_4$ where $\operatorname{gcd}(b_3,p_i)=\operatorname{gcd}(b_4,p_i)=1$ for all $i$ such that 
    \[ \mathrm{N}_{F/\Q}(\mathrm{Tr}_{K/F}(\alpha)) = \frac{\mathrm{N}_{F/\Q}(c_2\mathrm{Tr}_{K/F}(\alpha)\I_F)}{\mathrm{N}_{F/\Q}(c_2\I_F)}=\prod_{i=1}^t p_i^{(a_{i3}-a_{i4})f_{p_i}}b_3b_4^{-1}.\]

    The third term is $$|\Delta_F| = \prod_{i=1}^tp_i^{a_{i5}f_{p_i}}b_5$$ where $a_{i5}\geq 0$ and positive integer $b_5$ such that $\operatorname{gcd}(b_5,p)=1$. 

   This gives 
    \begin{equation*}
        \begin{split}
            \prod_{i=1}^tp_i^{c_i}&= \mathrm{N}_{F/\Q}(\mathrm{Tr}_{K/F}(\alpha))\mathrm{N}_{F/\Q}(I_F)^2|\Delta_F|\\
            &=(\prod_{i=1}^tp_i^{(a_{i3}-a_{i4})f_{p_i}}b_3b_4^{-1})(\prod_{i=1}^tp_i^{2(a_{i1}-a_{i2})f_{p_i}}b^2_1b^{-2}_2)(\prod_{i=1}^tp_i^{a_{i5}f_{p_i}}b_5)\\
           &=\prod_{i=1}^tp_i^{2(a_{i1}-a_{i2})f_{p_i}+(a_{i3}-a_{i4})f_{p_i}+a_{i5}f_{p_i}}b^2_1b^{-2}_2b_3b_4^{-1}b_5.
        \end{split}
    \end{equation*} and since $b_1,\ldots,b_5$ are all coprime to $p_i$ for each $i$, then
   \[\prod_{i=1}^tp_i^{c_i} = \prod_{i=1}^tp_i^{f_{p_i}(2(a_{i1}-a_{i2})+(a_{i3}-a_{i4})+a_{i5})} \implies f_{p_i}\vert c_i\] for all $i.$

Moreover, assume that $\mathrm{N}_{F/\Q}(\mathrm{Tr}_{K/F}(\alpha))$ is a square and $c_i$ are all positive odd integer. 
That $c_i$ are odd immediately implies that $f_{p_i}$ is odd for all $i$. The assumption that $\mathrm{N}_{F/\Q}(\mathrm{Tr}_{K/F}(\alpha))$ is a square implies that there exist $a'_i\in \Z$ and positive integers $b_3',b_4' \in \Z$ where $\mathrm{gcd}(b_3',p_i)=\mathrm{gcd}(b_4',p_i)=1$ for all $i$ such that 
\[\mathrm{N}_{F/\Q}(\mathrm{Tr}_{K/F}(\alpha)) = \prod_{i=1}^tp_i^{(a_{i3}-a_{i4})f_{p_i}}b_3b_4^{-1}= (\prod_{i=1}^tp_i^{a'_if_{p_i}}b_3'b_4'^{-1})^2.\] This gives 
\[ c_i= f_{p_i}(2(a_{i1}-a_{i2})+2a'_i+a_{i5}) \implies a_{i5} \geq 1 \text{ is odd for all $i.$}\]
 Then 
\[ b_1^2b_2^{-2}{b_3'}^2{b'_4}^{-2}b_5 = 1 \iff b_5 =b_1^{-2}b_2^{2}{b_3'}^{-2}{b_4'}^2 \iff b_5 \text{ is a square}.\] This says that $|\Delta_F|= b_5\prod_{i=1}^tp_i^{a_{i5}f_{p_i}}$ where $a_{i5}f_{p_i}$ are all odd and $b_5$ is a square that is coprime to $p_i$ for all $i$. 
\end{proof}
\begin{corollary}\label{cor: Necessary condition for subideal lattice to attain optimality}

    Let $F$ be a subfield of degree $k=2$ and $d_2(I,b_\alpha)=\prod_{i=1}^tp_i^{c_i}$ where $p_i$ are distinct primes, $c_i$ are all positive odd integer and $\mathrm{N}_{F/\Q}(\mathrm{Tr}_{K/F}(\alpha))$ is a square. 
    If  $\gamma^{ideal}_{n,2}(I,b_\alpha) = \gamma_{n,2}(I,b_\alpha)$ is attained at $F$, then $F = \Q\left(\sqrt{\pm \prod_{i=1}^t p_i}\right).$
\end{corollary}

\begin{proof}
    Since $k=2$, then $F=\Q(\sqrt{d})$ is a Galois extension where $d$ is squarefree. It is known that $\Delta_F = d$ if $d \equiv 1 \pmod 4$ and $\Delta_F = 4d$ if $d\equiv 2,3 \pmod 4$, and by Proposition \ref{Prop: Necessary condition for subideal lattice to attain optimality}, $|\Delta_F| = b\prod_{i=1}^tp_i^{a_if_{p_i}}$ where $a_if_{p_i}$ are all positive odd integers and $b$ is a square coprime to $p_i$ for all $i$. Since $f_{p_i}$ are all odd, $f_{p_i}\geq 1$ and $[F:\Q]=2$, then $f_{p_i}=1$ for all $i$ and hence, $|\Delta_F|=\prod_{i=1}^t p_i^{a_i}b.$

If $d \equiv 1 \pmod 4$, then $|\Delta_F| = |d| =\prod_{i=1}^t p_i$ as $d$ is squarefree. 
If $d \equiv 2,3 \pmod 4$, we distinguish two cases: 
\begin{itemize}
    \item 
    If there exists $i$ such that  $p_i = 2$, then $v_2(4d) = 2 + v_2(d).$ Since $d$ is squarefree, then $v_2(d)\in \{0,1\}$ which implies that $a_i = 2 +v_2(d) \in \{2,3\}$. Since $a_i$ is odd, then $2+v_2(d)=3 \iff v_2(d)= 1.$ Also, since $a_i$ are all odd, then $v_{p_i}(d)=1$ for all $i.$ This implies that $|d|= \prod_{i=1}^t p_i$
    \item 
     If $p_i \neq 2$ for all $i$, then $|d|= \prod_{i=1}^t p_i$ as $d$ is squarefree and so $b = 4.$ 
\end{itemize}
In either case, $F = \Q\left(\sqrt{\pm \prod_{i=1}^t p_i}\right).$
\end{proof}

Note that Proposition \ref{Prop: Necessary condition for subideal lattice to attain optimality} and Corollary \ref{cor: Necessary condition for subideal lattice to attain optimality} only hold for a Galois subfield $F$ of degree $k$ and does not require the larger field $K$ to be Galois. In general, if $\mathrm{Gal}(K/\Q)$ is not abelian, then the subfield need not be Galois over $\Q$. In our case, we will consider $K$ to be a cyclotomic or the maximal real subfield of a cyclotomic field which are both abelian extensions, hence all subfields are Galois over $\Q.$

To illustrate the above results, we consider the well known lattices $E_6$ and $E_8$.

\subsection{Some Examples of $k$-Subideal Lattices }
\label{Ex: E6 lattice}
Known lattices include $A_n, D_n, E_6, E_7, E_8$ and the Leech lattice $\Lambda_{24}$,  see \cite{conway2013sphere} for definitions and properties of these lattices, such as their canonical Gram matrices. We recall that when a Gram matrix is not canonical, but is instead obtained as the conjugation by an orthogonal matrix, we refer to its corresponding lattice as being rotated. 

 The following is an ideal lattice construction of the $E_6$ lattice \cite[p.77]{Bayer-Fluckiger:1999}. Set $K= \Q(\zeta_9)$ which is of degree 6. Choose the fractional ideal $I = \mathfrak{P}^{-4}$ where $\mathfrak{P}=(1-\zeta_9)$ is the unique prime ideal above $3$. Let $\alpha =1.$ The ideal lattice $(I,b_\alpha)$ is a rotated $E_6$ lattice.

We know that the lattice $E_6$ is even, that  $d_1(I,b_\alpha)=2$ and since $\gamma_{6,2}= 3^{\frac{2}{3}}$ is attained by $E_6$, then $d_2(I,b_\alpha)=3$. Finally, the value $d_3(I,b_\alpha)= 4$ can be found from \cite[Proposition 2]{sawatani2010note}. Below we provide another simple proof.
\begin{proposition}\label{Prop: minimum discriminant of rank 3 sublattice in E_6}
    We have  
    \[d_3(E_6) = 4\] and so,
    \[  \gamma_{6,3} \geq \gamma_{6,3}(E_6)=  \frac{2^2}{3^\frac{1}{2}} \approx 2.3094. \]
\end{proposition}
\begin{proof}
    Let $L$ be a sublattice of rank 3. By definition, 
    \[ \gamma_{3,1}\geq \frac{d_1(L)}{\operatorname{disc}(L)^{\frac{1}{3}}}\iff \operatorname{disc}(L)^{\frac{1}{3}}\geq d_1(L)\gamma_{3,1}^{-1}.\] Since $L$ is a sublattice of $E_6$, $d_1(L) \geq d_1(E_6)= 2$ and recalling that $\gamma_{3,1}=2^\frac{1}{3}$, we get 
    \[
    \operatorname{disc}(L)^{1/3}\geq 2 \cdot 2^{-1/3} \Rightarrow 
    \operatorname{disc}(L)\geq 4.\] There is a sublattice of rank 3 whose Gram matrix is  
    \[\begin{pmatrix}
         2&-1&0\\
        -1&2&-1\\
        0&-1&2\\
    \end{pmatrix}\]
   of discriminant 4. This gives the lower bound 
    \[ \gamma_{6,3}\geq \gamma_{6,3}(E_6) = \frac{2^2}{3^\frac{3}{6}} \approx 2.3094.\]
\end{proof}
\begin{example}
Let $(I,b_\alpha)$ be the $E_6$ lattice constructed above. For $k=1$, we know that $d_1(I,b_\alpha)=2$ and by Proposition \ref{Prop: Necessary condition for subideal to attain optimality rank 1}, $
\frac{Tr_{K/\Q}(\alpha)}{p^{v_p(Tr_{K/\Q}(\alpha))}}=
\frac{Tr_{K/\Q}(1)}{2}=\frac{6}{2^1}=3
$
is not a square, hence $\gamma^{ideal}_{6,1}(I,b_\alpha)\neq \gamma_{6,1}(I,b_\alpha).$ Instead, the free $\Z$-module $(\frac{\zeta^3}{3}+\frac{2}{3})\Z\subset I$ of rank 1 attains the optimal value $\gamma_{6,1} = \left(\frac{64}{3}\right)^{\frac{1}{6}}.$     

For $k=2$, since $\alpha=1$, then $\mathrm{N}_{F/\Q}(\mathrm{Tr}_{K/F}(\alpha))$ is a square for every subfield $F$ of degree 2.  By Corollary \ref{cor: Necessary condition for subideal lattice to attain optimality}, a necessary condition for $\gamma^{ideal}_{n,k}(I,b_\alpha)$ to obtain optimality at a quadratic subfield $F_1$ is to have $F_1 = \Q(\sqrt{\pm 3})$, which exists since $F_1 = \Q(y)= \Q[X]/(X^2+3X+9)\cong \Q(\sqrt{-3})$ is the unique subfield of $K$ of degree 2 with relative discriminant $\Delta_{K/F_{1}}= 27\I_{F_1}$. We know that the ideal below $\mathfrak{P}$ in $F_1$ is $\sqrt{-3}\I_F$ and $e_{\sqrt{-3}\I_K}= 3$, so by Proposition \ref{Prop: general factorization of ideal below}, $I_{F_1} = I \cap F_1 = ((\sqrt{-3})\I_{F_1})^{\left\lceil\frac{-4}{3} \right\rceil} = (\sqrt{-3}\I_{F_1})^{-1}$ and 

  \begin{equation*}
      \begin{split}
           \gamma^{ideal}_{6,2}(I,b_1)&= \frac{N_{F_1/\Q}(Tr_{K/F_1}(1))}{N_{K/\Q}(1)^\frac{2}{6}}N_{F_1/\Q}(27\I_{F_1})^{\frac{-1}{3}}3^{2(-1-\frac{2}{6}(-4))}\\
           &=3^{\frac{2}{3}}\\
           &= \gamma_{6,2}\\ &= \gamma_{6,2}(I,b_1).
      \end{split}
  \end{equation*}

  Set $F_2 =\Q(y)= \Q[X]/(X^3-3X+1)$ which is the unique subfield of $K_9$ of degree 3. By Proposition \ref{Prop: minimum discriminant of rank 3 sublattice in E_6}, we know that $d_3(I,b_\alpha) = 4 = 2^2$ and it can be checked that 2 is inert in $F_2$, hence $f_2=3$ which does not divide the exponent $2$. By Proposition \ref{Prop: Necessary condition for subideal lattice to attain optimality}, $\gamma^{ideal}_{6,3}(I,b_\alpha)\neq \gamma_{6,3}(I,b_\alpha)$, which we verify with $\mathrm{disc}(I_F)= 8 \neq 4.$ Instead, the free $\Z$-module generated by $\{\frac{1}{3}(y+1)\zeta_9,\frac{1}{3}(y+1),\frac{1}{3}y(y+1)\}$ of rank 3 attains the optimal value $\gamma_{6,3}(I,b_\alpha)=  \frac{2^2}{3^\frac{3}{6}}\approx 2.3094.$

\end{example}

The lattice $E_8$ is an exceptional lattice in dimension 8. It is the unique unimodular even lattice in this dimension. It can be constructed from number fields as follows. Set $K= \Q(\zeta_{24})$ of degree 8. Consider the integral ideal $\mathfrak{P}=(\zeta^6_{24} + \zeta^5_{24} -\zeta^2_{24} -\zeta_{24} - 1)\I_K$ which is a prime ideal in $\I_K$ and $\alpha= \frac{1}{12}$. Then $(\mathfrak{P},b_\alpha)$ is a unimodular even lattice, thus the $E_8$ lattice. 

\begin{example}
    Repeating a similar computation done in Example \ref{Ex: E6 lattice}, we see that there exist a $2$- and a $4$-subideal lattice that are optimal:

\begin{enumerate}

\item  $ \gamma^{ideal}_{8,2}(\mathfrak{P},b_\alpha)=\gamma_{8,2}(\mathfrak{P},b_\alpha)$: let $F_1 = \Q[x]/(X^2-4X+16) \cong \Q(\sqrt{- 3})$ and let $\mathfrak{p}$ be the unique prime below $\mathfrak{P}$ in $F_1$. Then $\mathfrak{P}_{F_1} = \mathfrak{p}$ and 
    \[\gamma_{8,2}((\mathfrak{P},b_\alpha)),(\mathfrak{P}_{F_1},b_\alpha|_{F_1}))= \operatorname{disc}(\mathfrak{P}_{F_1},b_\alpha|_{F_1})= 3= \gamma_{8,2}(\mathfrak{P},b_\alpha).\] This proves that $\gamma_{8,2}(\mathfrak{P},b_\alpha) = \gamma^{ideal}_{8,2}(\mathfrak{P},b_\alpha)$.

\item $\gamma_{8,4}^{ideal}(\mathfrak{P},b_\alpha)=\gamma_{8,4}(\mathfrak{P},b_\alpha)$: let $F_2 = \Q(w) = \Q[X]/(X^4-2X^2+4)$ and let $\mathfrak{p}$ be the unique prime below $\mathfrak{P}$ in $F_2.$ Then $\mathfrak{P}_{F_2}= \mathfrak{p}$ and 
    \[ \gamma_{8,4}((\mathfrak{P},b_\alpha),(\mathfrak{P}_{F_2},b_\alpha|_{F_2})) =\operatorname{disc}(\mathfrak{P}_{F_2},b_\alpha|_{F_2})=  4 = \gamma_{8,4}.\] This proves that $\gamma_{8,4}(\mathfrak{P},b_\alpha) = \gamma^{ideal}_{8,4}(\mathfrak{P},b_\alpha)$.

\end{enumerate} 
\end{example}

In the next example, we will use Corollary \ref{cor: Necessary condition for subideal lattice to attain optimality} to construct an ideal lattice where there exists a 2-subideal lattice that is optimal. For instance, if we would like to construct an ideal lattice $(I,b_\alpha)$ such that $d_2(I,b_\alpha)=3$ attained by a $2$-subideal lattice, then the number field $K$ must contain a quadratic subfield which is $\Q(\sqrt{\pm 3})$. Below, we fix the quadratic subfield $F =\Q(\sqrt{-3})$ and search for a cyclotomic field that contains $F$ as the unique quadratic subfield. 
\begin{example}\label{Example: zeta27}
    Set $K = \Q(\zeta_{27})$ which is of degree 18. It can be checked that the only quadratic subfield it contains is $F=\Q(\sqrt{-3})$. The different ideal of $K$ (see \cite[p.195]{neukirch2013algebraic} on the different ideal) is $\mathfrak{D}_K= \mathfrak{P}^{45}$ where $\mathfrak{P}$ is the unique prime above 3. Let $I=\mathfrak{P}^{-22}$ and consider the ideal lattice $(I,b_\alpha)$ where $\alpha =1$ of discriminant 3. Using Proposition \ref{Prop: general factorization of ideal below}, the ideal below is $I_F = (\sqrt{-3}\I_K)^{-2}$ and the discriminant of $(I_F,b_\alpha|_F)$ is 3. The ideal lattice $(I,b_\alpha)$ is an integral even lattice where $d_1(I,b_\alpha)=2$ which implies that $d_2(I,b_\alpha)= 3$ by Proposition \ref{Prop: general bounds in CM case}. Hence, $\gamma_{18,2}(I,b_\alpha)=\gamma^{ideal}_{18,2}(I,b_\alpha)= \frac{3}{3^\frac{2}{18}} \approx 2.65526$ and so, $\gamma_{18,2}\geq 2.65526.$
\end{example}

\section{Lower Bounds from Known Lattices}
\label{sec:lower}
In this section, we will provide lower bounds to $\gamma_{n,2}, \gamma_{9,3}$ and $\gamma_{9,4}$ where $n=5,7,9$ and $24$ using interesting known lattices. We present the algebraic constructions of the lattices used for dimension 5, 7 and 9, which are also of independent interest.

\subsection{The Case of $\gamma_{5,2}$}\label{subsec: dimension 5}

We present constructions of two 5-dimensional lattices arising from number fields. The $D_5$ lattice provides a natural benchmark for the lower bound of $\gamma_{5,2}$, since \[
1.7230 \approx \gamma_{5,2}(D_5)
\leq \gamma_{5,2}\leq 2,\]
as established in \cite{oggier2025rankin}. To identify other 5-dimensional lattices whose Rankin invariants may exceed $\gamma_{5,2}(D_5)$, we examine the list of interesting 5-dimensional lattices available on \cite{Catalogue}. Using their canonical Gram matrices, we compute the smallest discriminant among rank-2 sublattices by direct inspection which could possibly not be the minimum and subsequently, compute the Rankin invariant based on this discriminant. This comparison shows that $A_5^{+3}$ (see Definition \ref{Def: Coxeter lattice}) and $D_5^*$ are the only lattices in the list whose Rankin invariants could potentially exceed $\gamma_{5,2}(D_5)$. We will algebraically construct these two lattices and compute the respective actual minimum discriminant of a rank 2 sublattice.

Let $L\subseteq \R^n$ be a lattice. Define its dual \cite[Definition 1.3.1]{Martinet2010PerfectLI} to be 
\[ L^* = \{x\in L\otimes_{\Z}\R  \mid x\cdot y\in \Z \text{ for any } y \in L\}\subseteq \R^n.\] 

In the particular case where the lattice comes from a number field, that is we consider $\sigma_\alpha(M)$ (resp $(M,b_\alpha)$) where $M$ is a free $\Z$-submodule (resp fractional ideal) of $K$, 
define its trace dual to be 
\[M^{\sharp} = \{x \in M\otimes_{\Z}\R  \mid b_\alpha(x,y) \in \Z \text{ for all } y\in M\}\subseteq K\otimes_\Q \R.\]

We note that $M\otimes_\Z\R \subseteq K\otimes_\Q \R$ is an $\R$-vector space of dimension at most $[K:\Q]$ where $K\otimes_\Q \R \cong \R^{[K:\Q]}$ as an $\R$-vector space. Using the canonical embedding $\sigma_\alpha$,  we identify any $m \otimes r \in K \otimes_\Q \R$ as $r\sigma_\alpha(m)\in \R^{[K:\Q]}$ which can be extended linearly to any element in $K\otimes_\Q \R$ through the $\R$-vector space isomorphism. This allows us to identify $M^\sharp$ in $\R^{[K:\Q]}$: for any $x = \sum_{i=1}^t m_i\otimes r_i \in M^\sharp,$ we write $\sigma_\alpha(x)= \sum_{i=1}^t r_i\sigma_\alpha(m_i) \in \R^{[K:\Q]}$, so $\sigma_\alpha(M^\sharp)$ is a lattice in $\R^{[K:\Q]}$.
\begin{lemma}\label{Lemma: trace dual equals to dual}
    For any free $\Z$-module $M$ in $K$, 
    \[\sigma_\alpha(M)^*= \sigma_\alpha(M^\sharp).\]
\end{lemma}
\begin{proof}
    By identifying $\sum_{i=1}^t x_i\otimes r_i \in M^\sharp \subset M\otimes_\Z \R$ as $\sum_{i=1}^t r_i\sigma_\alpha(x_i)\in \R^n$, we will show that $\sum_{i=1}^t x_i \otimes r_i \in M^\sharp \iff \sum_{i=1}^t r_i\sigma_\alpha(x_i)\in \sigma_\alpha(M)^*.$ Recall that from Subsection \ref{subsec:prelim}, we can view $b_\alpha$ as a $\R$-bilinear form on $K_\R$ and for any $x,y \in K$, then $b_\alpha(x,y)= \sigma_\alpha(x)\cdot \sigma_\alpha(y)$. We see that
  \begin{equation*}
      \begin{split}
          \sum_{i=1}^t  x_i\otimes r_i \in M^\sharp &\iff b_\alpha(\sum_{i=1}^t x_i\otimes r_i,y\otimes 1) \in \Z \text{ for all } y \in M,\\
          & \iff \sum_{i=1}^t r_i b_\alpha(x_i\otimes  1, y\otimes 1) \text{ for all } y \in M,\\
          &\iff \sum_{i=1}^t r_i\sigma_\alpha(x)\cdot \sigma_\alpha(y) \in \Z \text{ for all } y \in M,\\
          &\iff \sum_{i=1}^t r_i\sigma_\alpha(x_i) \in \sigma_\alpha(M)^*.
      \end{split}
  \end{equation*}
\end{proof}

Let $K=\Q(\zeta_m+\zeta^{-1}_m)$ be the maximal real subfield of $\Q(\zeta_m)$ of degree $n= \frac{\varphi(m)}{2}$ with the following canonical $\Q$-basis $\{\theta_1,\theta_2,\ldots,\theta_{n}\}$, where $\theta_i= \zeta^i_m+\zeta^{-i}_m$ for $i=1,\ldots, n.$ 

When $m=p$ is a prime, then $n=\frac{p-1}{2}$. Throughout this subsection and \ref{subsec: dimension 9} which discuss the case of dimension 5 and 9 respectively, let $\alpha = 2-\theta_1\in K$ which is totally positive, $N$ be the generator matrix of $(\I_K,b_\alpha)$ with rows $\sigma_\alpha(\theta_i)$ for $i=1,\ldots,n$ and 
\[ T = \begin{psmallmatrix}
    1&1&\ldots&1\\
    0&1&\ldots&1\\
    \vdots&\vdots&\ldots&\vdots\\
    0&0&\ldots&1
\end{psmallmatrix}.\] By \cite[Proposition 1]{bayer2004new}, $\frac{1}{\sqrt{p}}(\I_K,b_\alpha)$ is the $\Z^n$ lattice and 
\[ \frac{1}{p}(TN)(TN)^\top = I_n\] where $I_n$ is the $n\times n $ identity matrix.

\begin{proposition}\label{Prop: A_5 and A^*_5 lattice}
    Set $K= \Q(\zeta_{13}+\zeta_{13}^{-1})$ which is a totally real number field of degree $6$, $\alpha = 2-\theta_1 \in K$ which is totally positive and let $M$ be the free $\Z$-module of rank $5$ spanned by $\{\theta_1,\ldots ,\theta_5\}\subseteq \I_{K}$, then $\frac{1}{\sqrt{13}}\sigma_\alpha(M)$ is a rotated $A_5$ lattice and $\sqrt{13}\sigma_\alpha(M^\sharp)$ is a rotated $A^*_5$ lattice where a $\Z$-basis of $M^\sharp$ is given by 
    \begin{equation*}
    \begin{split}
        e_1 &=\frac{5}{78}\theta_1  +\frac{2}{39} \theta_2 +\frac{1}{26}\theta_3 +\frac{1}{39}\theta_4+ \frac{1}{78}\theta_5 ,\\
        e_2& = \frac{2}{39}\theta_1 +\frac{4}{39}\theta_2+ \frac{1}{13}\theta_3+ \frac{2}{39}\theta_4+\frac{1}{39}\theta_5 ,\\
        e_3 &= \frac{1}{26}\theta_1 + \frac{1}{13}\theta_2 +\frac{3}{26}\theta_3 +\frac{1}{13}\theta_4 +\frac{1}{26}\theta_5 ,\\
        e_4 &=\frac{1}{39}\theta_1 +\frac{2}{39}\theta_2 +\frac{1}{13}\theta_3 +\frac{4}{39}\theta_4 +\frac{2}{39}\theta_5 ,\\
        e_5 &= \frac{1}{78}\theta_1 + \frac{1}{39}\theta_2+\frac{1}{26}\theta_3 + \frac{2}{39}\theta_4 +\frac{5}{78}\theta_5 .
    \end{split}
\end{equation*}
\end{proposition}\begin{proof}
Let $A_5$ be a 5-dimensional lattice with canonical generator matrix given by
\[ M_{A_5}= \begin{psmallmatrix}
        -1&1&0&0&0&0\\
        0&-1&1&0&0&0\\
        0&0&-1&1&0&0\\
        0&0&0&-1&1&0\\
        0&0&0&0&-1&1
    \end{psmallmatrix}.\]Consider the following matrix 
    \[ B = \frac{1}{\sqrt{13}}M_{A_5}TN.\] Then 
    \begin{equation*}
        \begin{split}
            B&=\frac{1}{\sqrt{13}}\begin{psmallmatrix}
        -1&&&&&\\
        &-1&&&&\\
        &&-1&&&\\
        &&&-1&&\\
        &&&&-1&0
    \end{psmallmatrix}\begin{psmallmatrix}
        \leftarrow \sigma_\alpha(\theta_1)\rightarrow\\
        \leftarrow \sigma_\alpha(\theta_2)\rightarrow\\
           \leftarrow \sigma_\alpha(\theta_3)\rightarrow\\
        \leftarrow \sigma_\alpha(\theta_4)\rightarrow\\
        \leftarrow \sigma_\alpha(\theta_5)\rightarrow\\
        \leftarrow \sigma_\alpha(\theta_6)\rightarrow\\
    \end{psmallmatrix}\\
    &= \frac{1}{\sqrt{13}}\begin{psmallmatrix}
        \leftarrow \sigma_\alpha(-\theta_1)\rightarrow\\
        \leftarrow \sigma_\alpha(-\theta_2)\rightarrow\\
           \leftarrow \sigma_\alpha(-\theta_3)\rightarrow\\
        \leftarrow \sigma_\alpha(-\theta_4)\rightarrow\\
        \leftarrow \sigma_\alpha(-\theta_5)\rightarrow\\
    \end{psmallmatrix}
        \end{split}
    \end{equation*} is a generator matrix for $\frac{1}{\sqrt{13}}\sigma_\alpha(M)$ such that
    \begin{equation*}
        \begin{split}
            BB^\top &=( \frac{1}{\sqrt{13}}M_{A_5}TN)( \frac{1}{\sqrt{13}}M_{A_5}TN)^\top\\
            &=M_{A_5}(\frac{1}{\sqrt{13}}TN)(\frac{1}{\sqrt{13}}TN)^\top M_{A_5}^\top\\
            &=M_{A_5}I_6M_{A_5}^\top\\
            &= G_{A_5}.
        \end{split}
    \end{equation*}  Hence, $\frac{1}{\sqrt{13}}\sigma_\alpha(M)$ is a rotated $A_5$-lattice.  

   Recall that for any lattice $cL\subset \R^n$ where $c$ is any positive real number, then $(cL)^*= \frac{1}{c}L^*.$ By Lemma \ref{Lemma: trace dual equals to dual},  $(\frac{1}{\sqrt{13}}\sigma_\alpha(M))^*= \sqrt{13}\sigma_\alpha(M^\sharp)$ is a rotated $A^*_5$ lattice. We can compute a $\Z$-basis for the free $\Z$-module $M^\sharp.$ Let $\{e_1,\ldots {e_5}\}$ be the $\Z$-basis for $M^\sharp$ which can be characterized by $b_\alpha( e_i, \theta_j)= \delta_{ij}$. For each $i=1,\ldots, 5$, since $e_i \in M\otimes_\Z\R$, which by definition implies that there exists $\lambda_{ik} \in \R$ such that $e_i = \sum_{k=1}^5 \theta_k\otimes \lambda_{ik}$. This gives 

   \[b_\alpha(e_i,\theta_j)= \sum_{k=1}^5\lambda_{ik}b_\alpha(\theta_k,\theta_j) = \delta_{ij} \iff \sum_{k=1}^5 \lambda_{ik}G_{kj}=\delta_{ij}\]
 where $G= (\mathrm{Tr}_{K/\Q}(\alpha \theta_i\theta_j)_{1\leq i,j\leq 5}$. This yields the following system $(\lambda_{ij})_{1\leq i,j\leq 5}G = I_5 \iff (\lambda_{ij})_{1\leq i,j\leq 5} = G^{-1}$ where 
\[ G^{-1} = \begin{psmallmatrix}
    \frac{5}{78}&\frac{2}{39}&\frac{1}{26}&\frac{1}{39}&\frac{1}{78}\\[2mm]
    \frac{2}{39}&\frac{4}{39}&\frac{1}{13}&\frac{2}{39}&\frac{1}{39}\\[2mm]
    \frac{1}{26}&\frac{1}{13}&\frac{3}{26}&\frac{1}{13}&\frac{1}{26}\\[2mm]
    \frac{1}{39}&\frac{2}{39}&\frac{1}{13}&\frac{4}{39}&\frac{2}{39}\\[2mm]
    \frac{1}{78}&\frac{1}{39}&\frac{1}{26}&\frac{2}{39}&\frac{5}{78}\\[2mm]
\end{psmallmatrix}.\]
Therefore,
 \begin{equation*}
    \begin{split}
        e_1 &=\frac{5}{78}\theta_1  +\frac{2}{39} \theta_2 +\frac{1}{26}\theta_3 +\frac{1}{39}\theta_4+ \frac{1}{78}\theta_5 ,\\
        e_2& = \frac{2}{39}\theta_1 +\frac{4}{39}\theta_2+ \frac{1}{13}\theta_3+ \frac{2}{39}\theta_4+\frac{1}{39}\theta_5 ,\\
        e_3 &= \frac{1}{26}\theta_1 + \frac{1}{13}\theta_2 +\frac{3}{26}\theta_3 +\frac{1}{13}\theta_4 +\frac{1}{26}\theta_5 ,\\
        e_4 &=\frac{1}{39}\theta_1 +\frac{2}{39}\theta_2 +\frac{1}{13}\theta_3 +\frac{4}{39}\theta_4 +\frac{2}{39}\theta_5 ,\\
        e_5 &= \frac{1}{78}\theta_1 + \frac{1}{39}\theta_2+\frac{1}{26}\theta_3 + \frac{2}{39}\theta_4 +\frac{5}{78}\theta_5 .
    \end{split}
\end{equation*}
\end{proof}

\begin{definition}\label{Def: Coxeter lattice}\cite[p.53]{conway1988low}
    The lattice $A_n^{+r}$ is the unique sublattice of $A^*_n$ containing $A_n$ to index $r$ where $r$ divides $n+1.$
\end{definition}

We will provide an algebraic construction of $A_5^{+3}$ and show that its Rankin invariant does not outperform the $D_5$ lattice.

\begin{proposition}\label{Prop: Coxeter lattice construction}
    Set $K = \Q(\zeta_{13}+\zeta_{13}^{-1})$ which is a totally real number field of degree $6$, $\alpha = 2-\theta_1\in K$ which is totally positive and let $L$ be the free $\Z$-module in $K$ spanned by $\{f_1,\ldots,f_4,2f_5\}$ where 
    \[ f_i = \frac{1}{\sqrt{13}}\theta_i \text{ for } i=1,\ldots,4 \text{ and } f_5 = \frac{1}{\sqrt{13}}(-\frac{1}{6}\theta_1-\frac{1}{3}\theta_2-\frac{1}{2}\theta_3-\frac{2}{3}\theta_4+\frac{1}{6}\theta_5).\] Then $\sigma_\alpha(L)$ is a rotated $A^{+3}_5$ lattice
\end{proposition}\begin{proof}
  Note that $\sqrt{13} \in K$. To construct $A^{+3}_5$, using Proposition \ref{Prop: A_5 and A^*_5 lattice}, our task is to find a free $\Z$-submodule $L\subset K$ of rank $5$ such that $\frac{1}{\sqrt{13}}M \subseteq L \subseteq \sqrt{13}M^\sharp$ where $|L/\frac{1}{\sqrt{13}}M|=3.$ Firstly, we know that 
     \[ \begin{psmallmatrix}
         e_1\\
         e_2\\
         e_3\\
         e_4\\
         e_5
     \end{psmallmatrix} = \begin{psmallmatrix}
          \frac{5}{78}&\frac{2}{39}&\frac{1}{26}&\frac{1}{39}&\frac{1}{78}\\[2mm]
    \frac{2}{39}&\frac{4}{39}&\frac{1}{13}&\frac{2}{39}&\frac{1}{39}\\[2mm]
    \frac{1}{26}&\frac{1}{13}&\frac{3}{26}&\frac{1}{13}&\frac{1}{26}\\[2mm]
    \frac{1}{39}&\frac{2}{39}&\frac{1}{13}&\frac{4}{39}&\frac{2}{39}\\[2mm]
    \frac{1}{78}&\frac{1}{39}&\frac{1}{26}&\frac{2}{39}&\frac{5}{78}\\[2mm]
     \end{psmallmatrix}\begin{psmallmatrix}
         \theta_1\\ \theta_2\\ \theta_3\\ \theta_4\\ \theta_5
     \end{psmallmatrix}\]if and only if 
     \[ \begin{psmallmatrix}
         \sqrt{13}e_1\\
          \sqrt{13}e_2\\
           \sqrt{13}e_3\\
            \sqrt{13}e_4\\
             \sqrt{13}e_5\\
     \end{psmallmatrix}= 13 \begin{psmallmatrix}
          \frac{5}{78}&\frac{2}{39}&\frac{1}{26}&\frac{1}{39}&\frac{1}{78}\\
    \frac{2}{39}&\frac{4}{39}&\frac{1}{13}&\frac{2}{39}&\frac{1}{39}\\
    \frac{1}{26}&\frac{1}{13}&\frac{3}{26}&\frac{1}{13}&\frac{1}{26}\\
    \frac{1}{39}&\frac{2}{39}&\frac{1}{13}&\frac{4}{39}&\frac{2}{39}\\
    \frac{1}{78}&\frac{1}{39}&\frac{1}{26}&\frac{2}{39}&\frac{5}{78}\\
     \end{psmallmatrix}\begin{psmallmatrix}
        \frac{1}{\sqrt{13}}\theta_1\\ \frac{1}{\sqrt{13}}\theta_2\\\frac{1}{\sqrt{13}}\theta_3\\\frac{1}{\sqrt{13}}\theta_4\\\frac{1}{\sqrt{13}}\theta_5
     \end{psmallmatrix},\]
     where 
     \[ A = \left(13\begin{psmallmatrix}
          \frac{5}{78}&\frac{2}{39}&\frac{1}{26}&\frac{1}{39}&\frac{1}{78}\\
    \frac{2}{39}&\frac{4}{39}&\frac{1}{13}&\frac{2}{39}&\frac{1}{39}\\
    \frac{1}{26}&\frac{1}{13}&\frac{3}{26}&\frac{1}{13}&\frac{1}{26}\\
    \frac{1}{39}&\frac{2}{39}&\frac{1}{13}&\frac{4}{39}&\frac{2}{39}\\
    \frac{1}{78}&\frac{1}{39}&\frac{1}{26}&\frac{2}{39}&\frac{5}{78}\\
     \end{psmallmatrix}\right)^{-1} =
         \begin{psmallmatrix}
 2 & -1 &  0 &  0 &  0 \\
-1 &  2 & -1 &  0 &  0 \\
 0 & -1 &  2 & -1 &  0 \\
 0 &  0 & -1 &  2 & -1 \\
 0 &  0 &  0 & -1 &  2
\end{psmallmatrix}.\]
    There exist unimodular $5\times 5$ matrices $P,Q$ such that 
 \[ PAQ = \operatorname{diag}(1,1,1,1,6) = D\] where $D$ is the Smith Normal Form of $A$ since we know the quotient $A^*_5/A_5$ is a cyclic group of order $6.$ As $Q$ is unimodular, then 
 \[ Q^{-1}\begin{psmallmatrix}
     \sqrt{13}e_1\\
      \sqrt{13}e_2\\
      \sqrt{13}e_3\\
     \sqrt{13} e_4\\
     \sqrt{13} e_5
 \end{psmallmatrix}= \begin{psmallmatrix}
     f_1\\
     f_2\\
     f_3\\
     f_4\\
     f_5
 \end{psmallmatrix}\]is still a $\Z$-basis for $\sqrt{13}M^\sharp.$ Similarly, since $P$ is unimodular, then 
\[ P\begin{psmallmatrix}
     \frac{1}{\sqrt{13}}\theta_1\\
     \frac{1}{\sqrt{13}}\theta_2\\
     \frac{1}{\sqrt{13}}\theta_3\\
     \frac{1}{\sqrt{13}}\theta_4\\
     \frac{1}{\sqrt{13}}\theta_5
 \end{psmallmatrix}= \begin{psmallmatrix}
     g_1\\
     g_2\\
     g_3\\
     g_4\\
     g_5
 \end{psmallmatrix}\] is still a $\Z$-basis for $\frac{1}{\sqrt{13}}M$. This gives 
 \[ \begin{psmallmatrix}
     \frac{1}{\sqrt{13}}\theta_1\\
     \frac{1}{\sqrt{13}}\theta_2\\
     \frac{1}{\sqrt{13}}\theta_3\\
    \frac{1}{\sqrt{13}} \theta_4\\
     \frac{1}{\sqrt{13}}\theta_5
 \end{psmallmatrix} = A \begin{psmallmatrix}
     \sqrt{13}e_1\\
     \sqrt{13}e_2\\
    \sqrt{13} e_3\\
     \sqrt{13}e_4\\
     \sqrt{13}e_5
 \end{psmallmatrix} \iff P\begin{psmallmatrix}
     \frac{1}{\sqrt{13}}\theta_1\\
     \frac{1}{\sqrt{13}}\theta_2\\
     \frac{1}{\sqrt{13}}\theta_3\\
     \frac{1}{\sqrt{13}}\theta_4\\
     \frac{1}{\sqrt{13}}\theta_5
 \end{psmallmatrix} = DQ^{-1}\begin{psmallmatrix}
     \sqrt{13}e_1\\
     \sqrt{13}e_2\\
    \sqrt{13} e_3\\
     \sqrt{13}e_4\\
     \sqrt{13}e_5
 \end{psmallmatrix}  \iff \begin{psmallmatrix}
      g_1\\
     g_2\\
     g_3\\
     g_4\\
     g_5
 \end{psmallmatrix}= D \begin{psmallmatrix}
     f_1\\
     f_2\\
     f_3\\
     f_4\\
     f_5
 \end{psmallmatrix}.\] Hence, a $\Z$-basis for $\sqrt{13}M^\sharp$ and $\frac{1}{\sqrt{13}}M$ is $\{f_1,f_2,f_3,f_4,f_5\}$ and $\{f_1,f_2,f_3,f_4,6f_5\}$ respectively. To compute $Q$, we see that
 \begin{equation*}
     \begin{split}
         \begin{psmallmatrix}
     f_1\\
     f_2\\
     f_3\\
     f_4\\
     f_5
\end{psmallmatrix} &= D^{-1} \begin{psmallmatrix}
      g_1\\
     g_2\\
     g_3\\
     g_4\\
     g_5
 \end{psmallmatrix}\\
 &=Q^{-1}A^{-1}P^{-1}\begin{psmallmatrix}
      g_1\\
     g_2\\
     g_3\\
     g_4\\
     g_5
 \end{psmallmatrix}\\
 &= Q^{-1}A^{-1} 
\begin{psmallmatrix}
    \frac{1}{\sqrt{13}}\theta_1\\
     \frac{1}{\sqrt{13}}\theta_2\\
     \frac{1}{\sqrt{13}}\theta_3\\
    \frac{1}{\sqrt{13}} \theta_4\\
     \frac{1}{\sqrt{13}}\theta_5
\end{psmallmatrix}\\
&= \begin{psmallmatrix}
    1&0&0&0&0\\
    0&1&0&0&0\\
    0&0&1&0&0\\
    0&0&0&1&0\\
    -\frac{1}{6}&-\frac{1}{3}&-\frac{1}{2}&-\frac{2}{3}&\frac{1}{6}
\end{psmallmatrix}\begin{psmallmatrix}
    \frac{1}{\sqrt{13}}\theta_1\\
     \frac{1}{\sqrt{13}}\theta_2\\
     \frac{1}{\sqrt{13}}\theta_3\\
    \frac{1}{\sqrt{13}} \theta_4\\
     \frac{1}{\sqrt{13}}\theta_5
\end{psmallmatrix}.
     \end{split}
 \end{equation*}
 This gives
\begin{equation*}
    \begin{split}
        f_1 &= \frac{1}{\sqrt{13}}\theta_1,~
        f_2 =\frac{1}{\sqrt{13}} \theta_2,~
        f_3 = \frac{1}{\sqrt{13}}\theta_3,~
        f_4 = \frac{1}{\sqrt{13}}\theta_4\\
        f_5 &= \frac{1}{\sqrt{13}}(-\frac{1}{6}\theta_1-\frac{1}{3}\theta_2-\frac{1}{2}\theta_3-\frac{2}{3}\theta_4+\frac{1}{6}\theta_5).
    \end{split}
\end{equation*} Let $L$ be the free $\Z$-module spanned by $\{f_1,\ldots,2f_5\}$, then a Gram matrix for $\sigma_\alpha(L)$ is

\[ \begin{psmallmatrix}
    2&-1&0&0&0\\
    -1&2&-1&0&0\\
    0&-1&2&-1&0\\
    0&0&-1&2&-2\\
    0&0&0&-2&\frac{10}{3}
\end{psmallmatrix}\] with determinant $\frac{2}{3}$. This gives $\sqrt{\frac{\operatorname{disc}(\frac{1}{\sqrt{13}}\sigma_\alpha(M))}{\operatorname{disc}(\sigma_\alpha(L))}}= 3$ where $\frac{1}{\sqrt{13}}M\subset L$. Since $A^{+3}_5$ is the unique sublattice of $A^*_5$ containing $A_5$ of index 3, then $\sigma_\alpha(L)$ is a rotated $A^{+3}_5$ lattice.
\end{proof}

\begin{proposition}

    \[ d_2(A^{+3}_5) = \frac{4}{3}.\]
\end{proposition}\begin{proof}
Firstly, we normalize $\sigma_\alpha(L)$ to make it an integral lattice. It is known that $\sqrt{3}\sigma_\alpha(L)$ is an integral lattice \cite[Proposition 5.2.3]{Martinet2010PerfectLI}, which is also clear from its Gram matrix computed in Corollary \ref{Prop: Coxeter lattice construction} and $d_1(\sqrt{13}\sigma_\alpha(L))=d_1(\sqrt{3}A^{+3}_5)= 4$ \cite[Theorem 5.2.1]{Martinet2010PerfectLI}. By Corollary \ref{Cor: Possible values of d_2(L)}, we have that 
\[ d_2(\sqrt{3}\sigma_\alpha(L))\geq \frac{3}{4}d_1(\sqrt{3}\sigma_\alpha(L))= 12.\] 
Consider the free $\Z$-submodule of $\sqrt{3}L$ of rank 2 spanned by $\{x_1=\sqrt{3}(f_4+2f_5), x_2=\sqrt{3}(f_2+2f_3+3f_4+4f_5)\}$, then 
\[b_\alpha(x_1,x_1) = b_\alpha(x_2,x_2) = 4 \text{ and }b_\alpha(x_1,x_2)= 2\] and the discriminant of this sublattice is  
\[ (4)^2 - (2)^2 = 12.\] This gives
\[d_2(\sqrt{3}\sigma_\alpha(L)) \leq 12.\]Hence, $d_2(\sigma_\alpha(L)) = \frac{12}{3^2}= \frac{4}{3}.$

\end{proof}

\begin{corollary}\label{Cor: Coxeter lattice perform worst than D_5}
\[1.5681 \approx \frac{\frac{4}{3}}{(2/3)^\frac{2}{5}}=\gamma_{5,2}(A^{+3}_5) < \gamma_{5,2}(D_5) = \frac{3}{4^{\frac{2}{5}}}\approx 1.7230.\]
\end{corollary}

By Corollary \ref{Cor: Coxeter lattice perform worst than D_5}, the lattice $A^{+3}_5$ does not outperform $D_5$ to approximate $\gamma_{5,2}.$

Next we will construct $D^*_5$. Firstly, it was proven that there is no fractional ideal in any number field $K$ that yields $D_5$ for any totally positive element $\alpha \in K$ \cite{jorge2013rotated}. So, the construction will come from a free $\Z$-module.
\begin{proposition}\label{Prop: D_5 and D^*_5 lattice}
     Set $K= \Q(\zeta_{11}+\zeta_{11}^{-1})$ which is a totally real number field of degree 5, $\alpha = 2-\theta_1 \in K$ which is totally positive and let $M$ be the $\Z$-module spanned by $\{m_1,\ldots,m_5\}$ where $m_1= -\theta_1-2\theta_2-2\theta_3-2\theta_4-2\theta_5$ and $m_i= \theta_{i-1}$ for $i=2,\ldots,5$. Then $\frac{1}{\sqrt{11}}\sigma_\alpha(M)$ is a rotated $D_5$ lattice and $\sqrt{11}\sigma_\alpha(M^\sharp)$ is a rotated $D_5^*$ lattice where a $\Z$-basis for $M^\sharp$ is given by
   \begin{equation*}
       \begin{split}
           e_1 &=\frac{5}{44} m_1 +\frac{3}{44}m_2 + \frac{3}{22}m_3  + \frac{1}{11}m_4+ \frac{1}{22} m_5 ,\\
          e_2 &= \frac{3}{44} m_1 + \frac{5}{44}m_2 + \frac{3}{22}m_3  + \frac{1}{11}m_4+ \frac{1}{22} m_5 ,\\
          e_3 &= \frac{3}{22}m_1  + \frac{3}{22}m_2 + \frac{3}{11}m_3  + \frac{2}{11}m_4+ \frac{1}{11} m_5 ,\\
          e_4 &=  \frac{1}{11}m_1+ \frac{1}{11}m_2 + \frac{2}{11}m_3  + \frac{2}{11}m_4+ \frac{1}{11}m_5 ,\\
          e_5 &=  \frac{1}{22} m_1 + \frac{1}{22}m_2 + \frac{1}{11}m_3  + \frac{1}{11}m_4+ \frac{1}{11}m_5 .
       \end{split}
   \end{equation*}
\end{proposition}
\begin{proof}
    A canonical generator matrix for the lattice $D_5$ is 
   \[ M_{D_5}= \begin{psmallmatrix}
       -1&-1&0&0&0\\
       1&-1&0&0&0\\
       0&1&-1&0&0\\
       0&0&1&-1&0\\
       0&0&0&1&-1
   \end{psmallmatrix}\] and let $G_{D_5}= M_{D_5}M_{D_5}^\top$. Another generator matrix for the rotated $D_5$ lattice is 
   \[\frac{1}{\sqrt{11}}M_{D_5}TN\] since 
   \[ \frac{1}{\sqrt{11}}M_{D_5}TN(\frac{1}{\sqrt{11}}M_{D_5}TN)^\top = G_{D_5}.\] So
   \begin{equation*}
       \begin{split}
           \frac{1}{\sqrt{11}}M_{D_5}TN &=\frac{1}{\sqrt{11}}\begin{psmallmatrix}
               -1&-2&-2&-2&-2\\
               1&0&0&0&0\\
               0&1&0&0&0\\
               0&0&1&0&0\\
               0&0&0&1&0\\
           \end{psmallmatrix}\begin{psmallmatrix}
               \leftarrow \sigma_\alpha(e_1)\rightarrow\\
               \leftarrow \sigma_\alpha(e_2)\rightarrow\\
               \leftarrow \sigma_\alpha(e_3)\rightarrow\\
               \leftarrow \sigma_\alpha(e_4)\rightarrow\\
               \leftarrow \sigma_\alpha(e_5)\rightarrow\\
           \end{psmallmatrix}\\
           &= \frac{1}{\sqrt{11}}\begin{psmallmatrix}
               \leftarrow\sigma_\alpha(-\theta_1-2\theta_2-\ldots-2\theta_5)\rightarrow\\
               \leftarrow \sigma_\alpha(\theta_1)\rightarrow \\
                \leftarrow \sigma_\alpha(\theta_2)\rightarrow \\
                 \leftarrow \sigma_\alpha(\theta_3)\rightarrow \\
                  \leftarrow \sigma_\alpha(\theta_4)\rightarrow \\
           \end{psmallmatrix}.
       \end{split}
   \end{equation*} Let $M$ be the free $\Z$-module spanned by the $\Z$-basis $\{m_1,\ldots,m_5\}$ where 
   \[m_1=-\theta_1-2\theta_2-\ldots-2\theta_5,m_2=\theta_1,m_3=\theta_2,m_4=\theta_3,m_5=\theta_4,\] then $\frac{1}{\sqrt{11}}\sigma_\alpha(M)$ is a rotated $D_5$ lattice.

   By Lemma \ref{Lemma: trace dual equals to dual}, $\frac{1}{\sqrt{11}}\sigma_\alpha(M)^*=\sqrt{11}\sigma_\alpha(M^{\sharp})$ is a rotated $D_5^*$ lattice. Let $\{e_1,\ldots,e_5\}$ be a basis for $M^{\sharp}$ characterised by $b_\alpha(e_i,m_j) =\delta_{ij}$ where $m_j$ denote the $\Z$-basis of $M$ stated above. For each $i$, let $e_i = \sum_{k=1}^5 m_k\otimes 
   \lambda_{ik},$ then 
   \[ b_\alpha(e_i,m_j)=\delta_{ij} \iff \sum_{k=1}^5\lambda_{ik}b_\alpha(m_k,m_j)=\delta_{ij}\iff \sum_{k=1}^5\lambda_{ik}G_{kj}=\delta_{ij},\] where $G = (Tr_{K/\Q}(\alpha m_im_j))_{1\leq i,j\leq 5}$ is the Gram matrix of $\sigma_\alpha(M)$, a rotated $D_5$ lattice scaled by $\sqrt{11}$. This gives $(\lambda_{ij})_{1\leq i,j\leq 5} = G^{-1}$ where 
   \[G^{-1}= \begin{pmatrix}
\frac{5}{44} & \frac{3}{44} & \frac{3}{22} & \frac{1}{11} & \frac{1}{22} \\[2pt]
\frac{3}{44} & \frac{5}{44} & \frac{3}{22} & \frac{1}{11} & \frac{1}{22} \\[2pt]
\frac{3}{22} & \frac{3}{22} & \frac{3}{11} & \frac{2}{11} & \frac{1}{11} \\[2pt]
\frac{1}{11} & \frac{1}{11} & \frac{2}{11} & \frac{2}{11} & \frac{1}{11} \\[2pt]
\frac{1}{22} & \frac{1}{22} & \frac{1}{11} & \frac{1}{11} & \frac{1}{11}
\end{pmatrix}.\] A $\Z$-basis for $M^\sharp$ is given by 
   \begin{equation*}
       \begin{split}
           e_1 &=\frac{5}{44} m_1 +\frac{3}{44}m_2 + \frac{3}{22}m_3  + \frac{1}{11}m_4+ \frac{1}{22} m_5 ,\\
          e_2 &= \frac{3}{44} m_1 + \frac{5}{44}m_2 + \frac{3}{22}m_3  + \frac{1}{11}m_4+ \frac{1}{22} m_5 ,\\
          e_3 &= \frac{3}{22}m_1  + \frac{3}{22}m_2 + \frac{3}{11}m_3  + \frac{2}{11}m_4+ \frac{1}{11} m_5 ,\\
          e_4 &=  \frac{1}{11}m_1+ \frac{1}{11}m_2 + \frac{2}{11}m_3  + \frac{2}{11}m_4+ \frac{1}{11}m_5 ,\\
          e_5 &=  \frac{1}{22} m_1 + \frac{1}{22}m_2 + \frac{1}{11}m_3  + \frac{1}{11}m_4+ \frac{1}{11}m_5 .
       \end{split}
   \end{equation*} Using this $\Z$-basis, a Gram matrix for $\sqrt{11}\sigma_\alpha(M^\sharp)$ is 
   \[
\begin{psmallmatrix}
\frac{5}{4} & \frac{3}{4} & \frac{3}{2} & 1 & \frac{1}{2} \\
\frac{3}{4} & \frac{5}{4} & \frac{3}{2} & 1 & \frac{1}{2} \\
\frac{3}{2} & \frac{3}{2} & 3 & 2 & 1 \\
1 & 1 & 2 & 2 & 1 \\
\frac{1}{2} & \frac{1}{2} & 1 & 1 & 1
\end{psmallmatrix}\] of determinant $\frac{1}{4}.$
\end{proof}

\begin{corollary}\label{Cor: minimum discriminant of dual of D_5}
    \[d_2(D^*_5) =1.\]
\end{corollary}
\begin{proof}
    From Proposition \ref{Prop: D_5 and D^*_5 lattice}, we can normalize the Gram matrix found for $\sqrt{11}\sigma_\alpha(M^\sharp)$ by considering $2\sqrt{11}\sigma_\alpha(M^\sharp)$ which gives the following integral Gram matrix
    \[\begin{psmallmatrix}
        5&3&6&4&2\\
        3&5&6&4&2\\
        6&6&12&8&4\\
        4&4&8&8&4\\
        2&2&4&4&4
    \end{psmallmatrix}.\] Since $d_1(2D^*_5)=4$ and there exist two linearly independent lattice points of squared norm $4$, by Corollary \ref{Cor: Possible values of d_2(L)},
    \[ 12=\frac{3}{4}(4)^2\leq d_2(2\sqrt{11}\sigma_\alpha(M^\sharp)) \leq 4^2= 16.\] By direct inspection of the Gram matrix, we see that 
    \[ d_2(2\sqrt{11}(M^\sharp,b_\alpha)) \leq  \mathrm{det}\begin{pmatrix}
        5&3\\3&5
    \end{pmatrix}= 16.\] Let $L\subset 2\sqrt{11}M^\sharp$ be a free $\Z$-submodule of rank 2. By Lemma \ref{Lemma: Lemma for modifying basis}, there exists a $\Z$-basis $\{x,y\}\subset L$ for $L$ where we denote $\ell_1 = 44b_\alpha(x,x), \ell_2=44b_\alpha(y,y)$ and $b=44b_\alpha(x,y)$ such that 
    \[\ell_1\leq \ell_2 \text{ and } |b|\leq \frac{\ell_1}{2}.\] In the case where $\ell_1=4$ and $\ell_2 \geq 5$, then 
    \[ \operatorname{disc}(L)= \ell_1\ell_2 -b^2 \geq 4(5)-2^2 =16.\] In the case where $\ell_1 ,\ell_5 \geq 5$, then 
    \[ \operatorname{disc}(L)\geq \frac{3}{4}\ell_1^2 \geq 18.75>16.\] Lastly, we will show that if $\ell_1=\ell_2= 4$, then $\mathrm{disc}(L)\geq 16$. An alternative definition \cite[p.4]{elkies2009theta} of the $D^*_5$ lattice is $\Z^5 \cup ((\frac{1}{2},\frac{1}{2},\frac{1}{2},\frac{1}{2},\frac{1}{2})+\Z^5)$. So, $2D^*_5$ lattice can be described as $2\Z^5\cup((1,1,1,1,1)+2\Z^5)$. Clearly, if $x \in 2\Z^5$, then $x = (2x_1,\ldots,2x_5)$ and if $\|x\|^2 =4$, then 
    \[ 4= 4\sum_{i=1}^5x_i^2\iff \sum_{i=1}^5x^2_i = 1\] which is attained exactly by the vectors $\pm e_i$ where $e_i$ takes the value $1$ at the $i$-th position and 0 everywhere. This contributes to $10$ lattice points of squared norm $4$. If $x \in ((1,1,1,1,1)+2\Z^5)$ and $\|x\|^2 = 4$, then 
    \[ 4 = \sum_{i=1}^5(2x_i+1)^2\geq 5\] since $(2x_i+1)^2 \geq 1$ as $2x_i+1$ is odd for all $i=1,\ldots,5$, hence a contradiction. This means that there are no minimum vectors contributed from $(1,1,1,1,1)+2\Z^5$. We know that for all $i\neq j$, $\pm e_i\cdot \pm e_j =0$. Since $2\sqrt{11}\sigma_\alpha(M^\sharp)$ is the rotated $2D^*_5$ lattice, then
    \[ \operatorname{disc}(L) = 4^2  -0 = 16.\] In all cases, 
    \[ \operatorname{disc}(L)\geq 16\] which implies that 
    \[ d_2(D^*_5)= d_2(\sqrt{11}(M^\sharp,b_\alpha)) = 1.\] 
     
\end{proof}

\begin{corollary}\label{Cor: New approximation for dimension 5}
    \[1.7411 \approx\frac{1}{(\frac{1}{4})^\frac{2}{5}}\leq \gamma_{5,2} \leq 2.\]
\end{corollary}
\begin{proof}
    Combining with the upper bound from \cite[Corollary 5]{oggier2025rankin}, we have 
    \[ 1.7411 \approx\frac{1}{(\frac{1}{4})^\frac{2}{5}}=\gamma_{5,2}(D_5^*)\leq \gamma_{5,2}\leq 2\]
\end{proof}

Corollary \ref{Cor: New approximation for dimension 5} improves the known lower bound found in $\gamma_{5,2}$ \cite[Corollary 5]{oggier2025rankin}.

\subsection{The Case of $\gamma_{7,2}$}
An approximation of $\gamma_{7,2}$ using Proposition \ref{Prop: Approximating using lower rank lattice} is
\[\gamma_{7,2} \geq\frac{3}{4}\gamma_{6,1}^\frac{12}{7}\approx 1.798.\] However, there is a better benchmark using the $D_7$ lattice \cite[Corollary 5]{oggier2025rankin}, which is
\[ 2.0189\approx \gamma_{7,2}(D_7) \leq \gamma_{7,2}\leq 3.1748. \]  In this subsection, we will show that the $E_7$ lattice is a better approximation for $\gamma_{7,2}.$
The following is a known ideal lattice construction of the $E_8$ lattice.
\begin{proposition}[\protect{\cite{watanabe2018e8}}]\label{Prop: E8 lattice}
        Set $K= \Q(\zeta_{24})$ which is a CM field of degree 8, $\alpha = \frac{1}{12}$ and $I=(1+\zeta_{24}^{8}+\zeta^{8}_{24}\zeta_{24}+\zeta^{8}_{24}\zeta^2_{24}).$ The ideal lattice $(I,b_\alpha)$ is even, integral and unimodular. Hence, isomorphic to the $E_8$ lattice.
\end{proposition}
A $\Z$-basis for $I$ is $\{e_0,\ldots,e_7\}$ where $e_0=3, e_1=\zeta_{24}^4+1, e_2=3\zeta_{24}, e_3 = \zeta^5_{24}+\zeta_{24}, e_4= \zeta^2_{24}+\zeta_{24}+1,e_5= \zeta^6_{24} + 2\zeta^2_{24} +1, \zeta^3_{24} + \zeta_{24} +1$ and $e_7 = \zeta^7_{24} +2\zeta_{24} +1$.

\begin{proposition}
    Let $(I,b_\alpha)$ be the $E_8$ lattice from Proposition \ref{Prop: E8 lattice} and $M$ be a free $\Z$-module with basis $\{e_0-3e_5, e_1-2e_5, e_2, e_3, e_4-2e_5, e_6-e_5, e_7-2e_5\}\subset I.$ Then $\sigma_\alpha(M)$ is a rotated $E_7$ lattice. Hence, 
    \[ 2.461=\frac{3}{2^{\frac{2}{7}}}\leq \gamma_{7,2} \leq 3.1748.\]
\end{proposition}
\begin{proof}
Let $ v \in (I,b_\alpha)$ such that $b_\alpha(v,v)= d_1(I,b_\alpha)=2$. Using the definition of $E_7$ \cite[p.124]{conway2013sphere}, then a rotated $E_7$ lattice is the following set
\[M = \{x \in I\mid b_{\alpha}(x,v)= 0 \}.\] A lattice point in $M$ of minimum squared norm is given for $v= e_1.$ Let $ x \in M$, then $x = \sum_{i=0}^7 a_i e_i$ where $a_i \in \Z$ for all $i.$ We compute that $b_\alpha(e_i,e_1)=3,2,0,0,2,1,1,2$ respectively. Hence, 
\[ b_\alpha(x,e_1)= \sum_{i=0}^7 a_iTr_{K/\Q}(\alpha e_i \overline{e_1}) = 0\] if and only if 
\[ 3a_0+2a_1+2a_4+a_5+a_6+2a_7 = 0.\] By letting 
\[ a_5 =-(3a_0+2a_1+2a_4+a_6+2a_7),\]  we see that 
\[ x=  a_0(e_0-3e_5)+1_1(e_1-2e_5)+a_2e_2 + a_3e_3 +a_4(e_4-2e_5)+a_6(e_6-e_5)+a_7(e_7-2e_5).\] Hence, a $\Z$-basis for $M$ is the set 
\[ \{e_0-3e_5, e_1-2e_5, e_2, e_3, e_4-2e_5, e_6-e_5, e_7-2e_5\}\subset I\] and $\sigma_\alpha(M)$ is a rotated $E_7$ lattice. 
\end{proof}

\begin{proposition}
    \[d_2(E_7)= 3.\]
\end{proposition}
\begin{proof}
     We know that $d_1(\sigma_\alpha(M))= d_1(E_7) =2$ and by Corollary \ref{Cor: Possible values of d_2(L)}, then 
\[ d_2(\sigma_\alpha)\geq \frac{3}{4}(2)^2 = 3.\] Notice that $b_\alpha(e_3,e_3)= 2$. Since the $24$-th cyclotomic polynomial is $X^8-X^4+1$, we have that 
\[e_3 \zeta^4_{24}= (\zeta_{24}^5 + \zeta_{24})\zeta^4_{24}= (\zeta_{24}^8 + \zeta_{24}^4)\zeta_{24} =(2\zeta_{24}^4-1)\zeta_{24} \] where $\zeta^4_{24}$ is a primitive $6$-th root of unity. We also have 
\[ -e_2 + 2e_3 = -3\zeta_{24}+2\zeta_{24}^{5}+2\zeta_{24}= (2\zeta_{24}^4-1)\zeta_{24}\] and hence, $e_3\zeta_{24}^4 \in M.$ By the proof of Corollary \ref{Prop: d_2 when 6th root of unity is present}, we know that $\{e_3,e_3\zeta_{24}\}$ is linearly independent such that $b_\alpha(e_3,e_3)=b_\alpha(e_3\zeta_{24}^4,e_3\zeta_{24}^4)= 2$ and $b_\alpha(e_3,e_3\zeta_{24}^4) =1$ which gives a rank 2 sublattice of discriminant $3,$ attaining the lower bound. Therefore \[d_2(E_7)= d_2(\sigma_\alpha(M))=3.\]
\end{proof}
\begin{corollary}
    \[ 2.4610\approx\frac{3}{2^{\frac{2}{7}}}\leq \gamma_{7,2}\leq 2^\frac{5}{3}\approx  3.1748.\]
\end{corollary}\begin{proof}
    Combining with the upper bound from \cite[Corollary 5]{oggier2025rankin}, 
    \[ 2.461\approx\frac{3}{2^{\frac{2}{7}}}=\gamma_{7,2}(E_7)\leq \gamma_{7,2} \leq 2^\frac{5}{3}\approx 3.1748.\]
\end{proof}

\subsection{The Cases of $\gamma_{9,2},\gamma_{9,3}$ and $\gamma_{9,4}$}\label{subsec: dimension 9}
The lattice $\Lambda_9$ can be defined using the $E_8$ lattice by \cite[Example 1]{sikiric2025lattice}
\[ \Lambda_9 = \{(x\vert 0)\in \R^9 \vert  x \in E_8\}+ \begin{psmallmatrix}
    1\\0\\0\\0\\0\\0\\0\\0\\1
\end{psmallmatrix}^\top\Z.\] Using this definition, we see that a generator matrix of $\Lambda_9$ is given by 
\[M_{\Lambda_9}=\begin{psmallmatrix}
        2&0&0&0&0&0&0&0&0\\
        -1&1&0&0&0&0&0&0&0\\
        0&-1&1&0&0&0&0&0&0\\
        0&0&-1&1&0&0&0&0&0\\
        0&0&0&-1&1&0&0&0&0\\
        0&0&0&0&-1&1&0&0&0\\
        0&0&0&0&0&-1&1&0&0\\
        \frac{1}{2}&\frac{1}{2}&\frac{1}{2}&\frac{1}{2}&\frac{1}{2}&\frac{1}{2}&\frac{1}{2}&\frac{1}{2}&0\\
        1&0&0&0&0&0&0&0&1
    \end{psmallmatrix}\] and let $G_{\Lambda_9}= M_{\Lambda_9}M_{\Lambda_9}^\top$. We note that $\Lambda_9$ is not an integral lattice, however, $\sqrt{2}\Lambda_9$ is.

    It was proven that the Hermite constant $\gamma_{9,1}=2$ is attained by the $\Lambda_9$ lattice \cite[Corollary 1]{sikiric2025lattice} and from Example \ref{Ex: dim 9 and 25}, we know that 
    \[ \gamma_{9,2}\geq 2.5717.\]Next, we provide an algebraic construction of the $\Lambda_9$ lattice and provide bounds for $\gamma_{9,2},\gamma_{9,3}$ and $\gamma_{9,4}.$
\begin{proposition}\label{Prop: laminated lattice}
    Set $K=\Q(\zeta_{19}+\zeta_{19}^{-1})$ which is a totally real number field of degree $9$, $\alpha = 2-\theta_1 \in K$ totally positive and let $M$ be the free $\Z$-module spanned by $\{e_1,\ldots,e_9\}\subset K$ where 
    \begin{equation*}
        \begin{split}
            e_1 &=2(\theta_1+\ldots+\theta_9),\\
            e_2 &= -\theta_1,\\
            e_3 &= -\theta_2,\\
            e_4 &= -\theta_3,\\
            e_5 &= -\theta_4,\\
            e_6 &= -\theta_5,\\
            e_7 &= -\theta_6,\\
            e_8 &=\frac{1}{2}\theta_1+\theta_2+\frac{3}{2}\theta_3+2\theta_4+\frac{5}{2}\theta_5 + 3\theta_6+\frac{7}{2}\theta_7+4\theta_8+4\theta_9,\\
            e_9 &= \theta_1+\ldots+\theta_8+2\theta_9.
        \end{split}
    \end{equation*} Then $\sqrt{\frac{1}{19}}\sigma_\alpha(M)$ is a rotated $\Lambda_9$ lattice.
\end{proposition}
\begin{proof}
   Since
    \[ M_{\Lambda_9}T = \begin{psmallmatrix}
        2&2&2&2&2&2&2&2&2\\
        -1&0&0&0&0&0&0&0&0\\
        0&-1&0&0&0&0&0&0&0\\
        0&0&-1&0&0&0&0&0&0\\
        0&0&0&-1&0&0&0&0&0\\
        0&0&0&0&-1&0&0&0&0\\
        0&0&0&0&0&-1&0&0&0\\
        0&0&0&0&0&0&-1&0&0\\
        \frac{1}{2}&1&\frac{3}{2}&2&\frac{5}{2}&3&\frac{7}{2}&4&4\\
        1&1&1&1&1&1&1&1&2
    \end{psmallmatrix}\] and 
    \[( \sqrt{\frac{1}{19}}M_{\Lambda_9}TN)(\sqrt{\frac{1}{19}}M_{\Lambda_9}TN)^\top= M_{\Lambda_9}M_{\Lambda_9}^\top= G_{\Lambda_9},\]
    then $\sqrt{\frac{1}{19}}M_{\Lambda_9}T N$ is a generator matrix of $\Lambda_9$. By direct computation, we also see that it is a generator matrix for $\sqrt{\frac{1}{19}}\sigma_\alpha(M)$ and hence, $\sqrt{\frac{1}{19}}\sigma_\alpha(M)$ is a rotated $\Lambda_9$ lattice. A Gram matrix of $\sqrt{\frac{1}{19}}\sigma_\alpha(M)$ is given by 
    \[\begin{psmallmatrix}
4 & -2 & 0 & 0 & 0 & 0 & 0 & 1 & 2\\
-2 & 2 & -1 & 0 & 0 & 0 & 0 & 0 & -1\\
0 & -1 & 2 & -1 & 0 & 0 & 0 & 0 & 0\\
0 & 0 & -1 & 2 & -1 & 0 & 0 & 0 & 0\\
0 & 0 & 0 & -1 & 2 & -1 & 0 & 0 & 0\\
0 & 0 & 0 & 0 & -1 & 2 & -1 & 0 & 0\\
0 & 0 & 0 & 0 & 0 & -1 & 2 & 0 & 0\\
1 & 0 & 0 & 0 & 0 & 0 & 0 & 2 & \frac{1}{2}\\
2 & -1 & 0 & 0 & 0 & 0 & 0 & \frac{1}{2} & 2
\end{psmallmatrix}.\] 
\end{proof}

\begin{proposition}\label{Prop: minimum of laminated lattice}
    \[d_2(\Lambda_9) = 3.\]
\end{proposition}
\begin{proof}
     We consider the integral lattice $\sqrt{\frac{2}{19}}\sigma_\alpha(M)$. Since $d_1(\sqrt{\frac{2}{19}}\sigma_\alpha(M)) = 4,$ by Corollary \ref{Cor: Possible values of d_2(L)}, we have that
\[ \gamma_{9,2}(\sqrt{\frac{2}{19}}\sigma_\alpha(M))\geq \frac{3}{4}d_1(\sqrt{\frac{2}{19}}\sigma_\alpha(M))^2 =12 \implies d_2(\sqrt{\frac{2}{19}}\sigma_\alpha(M))\geq 12.\] From the Gram matrix of $\sqrt{\frac{2}{19}}\sigma_\alpha(M)$ gotten by scaling the Gram matrix given in Proposition \ref{Prop: laminated lattice} by 2, there is a sublattice with Gram matrix $\begin{psmallmatrix}
    4&-2\\-2&4
\end{psmallmatrix}$ of determinant 12 and hence, $d_{2}(\sqrt{\frac{2}{19}}\sigma_\alpha(M))= 12$ which implies that $d_2(\sqrt{\frac{1}{19}}\sigma_\alpha(M))=3$. As $\sqrt{\frac{1}{19}}\sigma_\alpha(M)$ is a rotated $\Lambda_9$ lattice, then $d_2(\Lambda_9)=3.$
\end{proof}
From Proposition \ref{Prop: minimum of laminated lattice}, $\gamma_{9,2}(\Lambda_9)=\frac{3}{1^\frac{2}{9}}= 3$ which outperforms the lower bound given in Example \ref{Ex: dim 9 and 25}.
\begin{proposition}
\begin{itemize} For dimension $9$, we have the following bounds:
    \item $3\leq \gamma_{9,2}\leq 3(2)^\frac{1}{4} \approx 3.5676$,
    \item $4\leq \gamma_{9,3}\leq 2^\frac{19}{8}\approx 5.1874$,
    \item $4\leq \gamma_{9,4}\leq 2^\frac{5}{2} \approx 5.6569$.
\end{itemize}
    
\end{proposition}
\begin{proof}
  The  Rankin's inequality (see \cite[Theorem 2.8.6]{Martinet2010PerfectLI}) states that for any integers $k,h$ such that $0 \leq k \leq h \leq n$, then $\gamma_{n,k}\leq \gamma_{h,k}(\gamma_{n,h})^\frac{k}{h}$. Also, $\gamma_{n,k}=\gamma_{n,n-k}$ \cite[Proposition 2.8.5]{Martinet2010PerfectLI}.
\begin{itemize}
    \item We see that a lower bound can be attained using the $\Lambda_9$ lattice, so
\[ \gamma_{9,2} \geq \gamma_{9,2}(\Lambda_9)=3.\] Let $n= 9, h=8$ and $k =2$, we see that 
\[\gamma_{9,2}\leq \gamma_{8,2}(\gamma_{9,8})^\frac{2}{8}= \gamma_{8,2}(\gamma_{9,1})^\frac{1}{4}=3(2)^\frac{1}{4}\approx 3.5676\]
    \item  Using $n=9, h=3, k=1$, we see that 
  \[ \gamma_{9,3}\geq \left(\frac{\gamma_{9,1}}{\gamma_{3,1}}\right)^3= \left(\frac{2}{\sqrt[3]{2}}\right)^3 = 4\] and using $n=9, h=8, k=3$, we get 
  \[ \gamma_{9,3}\leq \gamma_{8,3}(\gamma_{9,8})^\frac{3}{8}= \gamma_{8,3}(\gamma_{9,1})^\frac{3}{8}=4(2)^\frac{3}{8}= 2^\frac{19}{8}.\] This gives 
  \[ 4\leq \gamma_{9,3}\leq 2^\frac{19}{8}\approx 5.1874.\]
  \item Lastly, for $n=9, h=4, k=1$ for the lower bound and $n=9, h=8, k=4$ for the upper bound, we get 
  \[ 4\leq \gamma_{9,4}\leq 2^\frac{5}{2}\approx 5.6569.\]
\end{itemize}
 
\end{proof}

\subsection{The Case of $\gamma_{24,2}$}
To the best of our knowledge, the following result has not been formally shown in the literature, thus we include it for completeness.
\begin{proposition}
    Let $\Lambda_{24} \subset \R^{24}$ be the leech lattice. Then 
    \[ d_2(\Lambda_{24}) = 12. \] In particular,
    \[12\leq  \gamma_{24,2}\leq 16.\]
\end{proposition}
\begin{proof}
    We know that the Leech lattice $\Lambda_{24}$ is an even unimodular lattice for which $d_1(\Lambda)=4$ \cite[p.131]{conway2013sphere}. The Leech lattice $\Lambda_{24}$ can be constructed as an ideal lattice coming from the cyclotomic field $\Q(\zeta_{84})$ by \cite[Proposition 3.4]{Bayer-Fluckiger:1999}. Since $6$ divides $84$, there exists a $6$-th root of unity in $\Z[\zeta_{84}]$ and by Proposition \ref{Prop: d_2 when 6th root of unity is present}, 
    \[ d_2(\Lambda_{24}) = \frac{3}{4}d_1(\Lambda_{24})^2= 12.\] Since $\gamma_{24,1}=4$ which is attained by $\Lambda_{24},$ using \cite[Proposition 2.8.4]{Martinet2010PerfectLI} that says $\gamma_{n,k}\leq \gamma_{n,1}^k$ where $n=24$ and $k=2$, we have that 
    \[12= \gamma_{24,2}(\Lambda_{24})\leq  \gamma_{24,2}\leq 4^2=16.\]
\end{proof}
\section{Open Questions}

In Section \ref{sec:subideal}, we showed that the bound $(n/k)^k$ of Proposition \ref{prop:nk} is achievable for $n=10$ and $k=2$. What are other values of $n,k$ for which the bound is achievable? Is it achievable for all values of $n,k$? Also, are there other families of ideal lattices coming from proper ideals where the ideal Rankin invariant coincide with the Rankin invariant?

In Section \ref{sec:loose}, we proved necessary condition for a $k$-subideal to attain optimality. In particular, for the case of $k=2$, as seen in Example \ref{Example: zeta27}, it is possible to construct an ideal lattice $(I,b_\alpha)$ such that $d_2(I,b_\alpha)=3$ by choosing a cyclotomic field containing the quadratic subfield $\Q(\sqrt{\pm 3})$. Since the Rankin constant is unknown for most dimensions, is it possible to utilize this necessary condition to construct ideal lattices in other dimensions with certain prescribed minimum discriminant to approximate the Rankin constant?

In Section \ref{sec:lower}, we showed that the Rankin invariants of $D_5^*$ and $E_7$ respectively provide a better approximation than what is known for $\gamma_{5,2}$ and $\gamma_{7,2}$ respectively. Table \ref{tab:summary of the bounds on rankin constant} shows a summary on the new bounds of $\gamma_{n,k}.$

\begin{table}
\[
\begin{array}{c|c}
n&\text{new bounds on } \gamma_{n,k} \\
\hline
5&\frac{1}{(\frac{1}{4})^\frac{2}{5}}\leq \gamma_{5,2} \leq 2\\ \hline
7& \frac{3}{2^{\frac{2}{7}}}\leq \gamma_{7,2} \leq 2^\frac{5}{3} \\
\hline
\multirow{3}{*}{9}&  3\leq \gamma_{9,2}\leq 3(2)^\frac{1}{4} \\
& 4\leq \gamma_{9,3}\leq 2^\frac{19}{8}\\
&4\leq \gamma_{9,4}\leq 2^{5/2}\\
\hline
24& 12 \leq \gamma_{24,2} \leq 16
\end{array}\]
\caption{\label{tab:summary of the bounds on rankin constant}
Summary of the new bounds on $\gamma_{n,2}$
}
\end{table}

\begin{conjecture}\label{Conj: dimension 5}
  $\gamma_{5,2}= \gamma_{5,2}(D^*_5) = \frac{1}{(\frac{1}{4})^\frac{2}{5}} \approx 1.7411.$
     
\end{conjecture} 
Assuming Conjecture \ref{Conj: dimension 5}, by using the identity $\gamma_{n,k}^{n-2k}\leq \gamma_{n-k,k}^{n-k}$ \cite[Corollary 2.8.9]{Martinet2010PerfectLI}: 
\[ \gamma_{7,2}^{3} \leq \gamma_{5,2}^5 \iff 2.461 \leq \gamma_{7,2} \leq \left(\frac{1}{(\frac{1}{4})^\frac{2}{5}}\right)^\frac{5}{3} \approx 2.5198.\]
\begin{conjecture}\label{Conj: dimension 7}
$\gamma_{7,2}= \gamma_{7,2}(E_7) = \frac{3}{2^\frac{2}{7}} \approx 2.461.$
\end{conjecture}

\section*{Acknowledgement}
 The work of Yuan Xiang Chew was supported by a Nanyang President's  Graduate Scholarship, whose support is gratefully acknowledged.
\bibliographystyle{plain}
\bibliography{references.bib} 
\end{document}